\documentclass[12pt,reqno]{amsart}

\usepackage{setspace,  
hyperref, epigraph} 

\usepackage{breakurl}   

\usepackage{tikz}

\usetikzlibrary{calc}

\numberwithin{equation}{section}

\numberwithin{figure}{section} 

\DeclareMathOperator{\Q}{{\mathbb Q}}

\DeclareMathOperator{\R}{{\mathbb R}}

\DeclareMathOperator{\N}{{\mathbb N}}

\newcommand{\st}{\textbf{st}}

\newcommand{\hr} {{{}^{\mathfrak{h}}\hspace*{-2.3pt}\R}}

\newcommand\astr{{{}^\ast\hspace{-2.5pt}\R}}

\author[J. Bair]{Jacques Bair}\address{J. Bair, HEC-ULG, University of
Liege, 4000 Belgium}\email{J.Bair@ULiege.be}

\author[P. B\l aszczyk]{Piotr B\l{}aszczyk}\address{P. B\l{}aszczyk,
Institute of Mathematics, Pedagogical University of Cracow,
Poland}\email{pb@up.krakow.pl}

\author[R. Ely]{Robert Ely}\address{R. Ely, Department of Mathematics,
University of Idaho, Moscow, ID 83844 US}\email{ely@uidaho.edu}

\author[P. Heinig]{Peter Heinig} \address{P. Heinig}
\email{heinig@ma.tum.de}

\author[M. Katz]{Mikhail G. Katz}\address{M. Katz, Department of
Mathematics, Bar Ilan University, Ramat Gan 52900
Israel}\email{katzmik@macs.biu.ac.il}


\begin{document}


\thispagestyle{empty}


\title{Leibniz's well-founded fictions and their interpretations}

\begin{abstract}
Leibniz used the term \emph{fiction} in conjunction with
infinitesimals.  What kind of fictions they were exactly is a subject
of scholarly dispute.  The position of Bos and Mancosu contrasts with
that of Ishiguro and Arthur.  Leibniz's own views, expressed in his
published articles and correspondence, led Bos to distinguish between
two methods in Leibniz's work: (A) one exploiting classical
`exhaustion' arguments, and (B) one exploiting \emph{inassignable}
infinitesimals together with a law of continuity.

Of particular interest is evidence stemming from Leibniz's work
\emph{Nouveaux Essais sur l'Entendement Humain} as well as from his
correspondence with Arnauld, Bignon, Dagincourt, Des Bosses, and
Varignon.  A careful examination of the evidence leads us to the
opposite conclusion from Arthur's.

We analyze a hitherto unnoticed objection of Rolle's concerning the
lack of justification for extending axioms and operations in geometry
and analysis from the ordinary domain to that of infinitesimal
calculus, and reactions to it by Saurin and Leibniz.

A newly released 1705 manuscript by Leibniz (\emph{Puisque des
personnes\ldots}) currently in the process of digitalisation, sheds
light on the nature of Leibnizian inassignable infinitesimals.

In a pair of 1695 texts Leibniz made it clear that his incomparable
magnitudes violate Euclid's Definition V.4, a.k.a.\;the
Archi\-me\-dean property, corroborating the non-Archimedean construal
of the Leibnizian calculus.

\end{abstract}

\maketitle
\tableofcontents

\epigraph{J'appelle \emph{grandeurs incomparables} dont l'une
multipli\'ee par quelque nombre fini que ce soit, ne s\c cauroit
exceder l'autre.  --G.\;W.\;Leibniz

\medskip\noindent 
[C]es touts infinis, et leurs oppos\'es infiniment
petits, ne sont de mise que dans le \emph{calcul des g\'eom\`etres},
tout comme les racines imaginaires de l'alg\`ebre.
--\emph{Th\'eophile}}

\section{Introduction}

Figure\;130 in l'Hospital's book \emph{Analyse des infiniment
petits\ldots} illustrates his \emph{Article\;163}.  This item concerns
the application of what is known today as l'H\^opital's rule to the
geometric situation of two curves crossing at a point~$B$ on the
$x$-axis.  

\begin{figure}
\begin{center}
\begin{tikzpicture}
  %
label (in the final picture) has parenthesis around the respective letter.
  \node[] (number) at (30pt,-35pt){$\textsl{Fig.\;130.}$};
  \node[scale=0.5] (remark) at (40pt,-42.5pt){$\text{%
%
%
}$};  
  %
    %
    \node[circle,draw,fill,inner sep=0em,minimum
size=0.75mm,black,opacity=1.0,label=left:$A$] (A) at (0pt,0pt){};
    \node[circle,draw,fill,inner sep=0em,minimum
size=0.75mm,black,opacity=1.0,label=below:$C$] (C) at (0pt,-100pt){};
    \draw[line width=0.25pt] (A) -- (C);
    %
    %
    %
    \node[circle,draw,fill,inner sep=0em,minimum
size=0.75mm,black,opacity=1.0,label=above:$M$] (M) at (90pt,123.0pt){};
    \node[circle,draw,fill,inner sep=0em,minimum
size=0.75mm,black,opacity=1.0,label=above right:$N$] (N) at
(90pt,60pt){};
    \node[circle,draw,fill,inner sep=0em,minimum
size=0.75mm,black,opacity=1.0,label=above right:$P$] (P) at (90pt,0pt){};
    \node[circle,draw,fill,inner sep=0em,minimum
size=0.75mm,black,opacity=1.0,label=below:$O$] (O) at (90pt,-79.5pt){};
    \draw[line width=0.5pt] (M)--(O);
    %
    \node[circle,draw,fill,inner sep=0em,minimum
size=0.75mm,black,opacity=1.0,label=above:$D$] (D) at (200pt,168.5pt){};
    \node[circle,draw,fill,inner sep=0em,minimum
size=0.75mm,black,opacity=1.0,label=below:$B$] (B) at (199.5pt,0pt){};
    \draw[line width=0.5pt] (D)--(B);
    %
    %
    \node[circle,draw,fill,inner sep=0em,minimum
size=0.75mm,black,opacity=1.0,label=above:$d$] (d) at (225pt,175pt){};
size=0.75mm,black,opacity=1.0,label=above left:$\text{\huge{\textsl{g}}}$]
(g) at (225pt,27.5pt){};
    \node[circle,draw,fill,inner sep=0em,minimum
size=0.75mm,black,opacity=1.0,label=above left:$g$] (g) at
(225pt,27.5pt){};    
    \node[circle,draw,fill,inner sep=0em,minimum
size=0.75mm,black,opacity=1.0,label=above right:$b$] (b) at
(225pt,0pt){};
    \node[circle,draw,fill,inner sep=0em,minimum
size=0.75mm,black,opacity=1.0,label=below left:$f$] (f) at
(225pt,-31pt){};
    \draw[line width=0.5pt,dashed] (d)--(f);    
    %
    %
    \draw[line width=0.25pt] (A)--($(b)+(25pt,0pt)$);    
    \draw (A) .. controls (10pt,100pt) and (125pt,150pt) ..
($(d)+(10pt,2.5pt)$);
    \draw (A) parabola bend (N) ($(f)+(10pt,-15pt)$);
    \draw (C) parabola ($(g)+(10pt,12pt)$);
  \end{tikzpicture}    
\end{center}
\caption{L'Hospital's Figure 130}
\label{f999}
\end{figure}
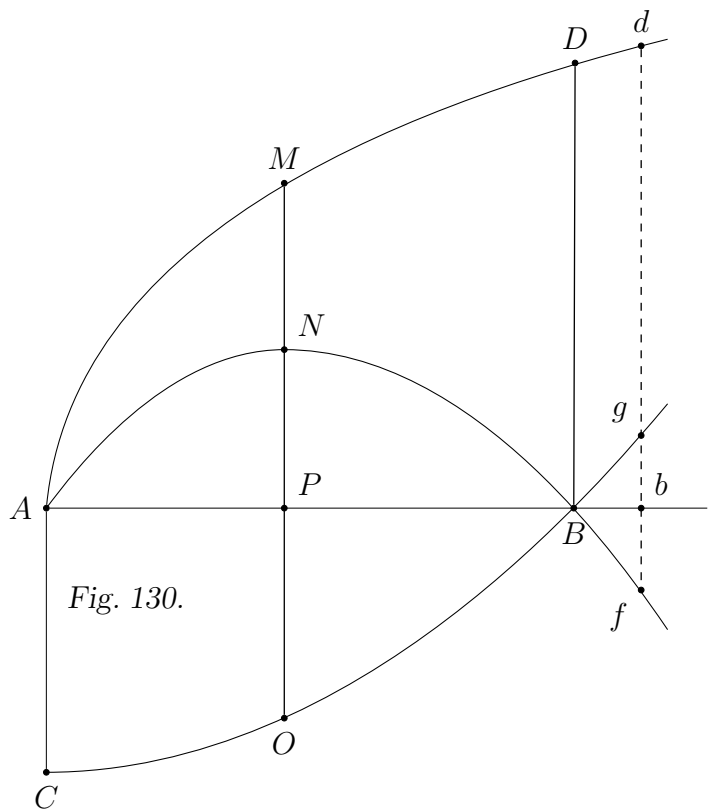

\subsection{Figure\;130} 
\label{s1}

L'Hospital's figure places another point~$b$ near~$B$ on the~$x$-axis,
as well as points~$f$ and~$g$ on the two curves with abscissa
(i.e.,~$x$-coordinate) equal to~$b$; see Figure~\ref{f999}.  The
author argues that~$bf$ and~$bg$ are the differentials corresponding
to the two curves.  He concludes that the ratio of the ordinates
(i.e.,~$y$-coordinates) of the two curves equals the ratio of their
differentials.

In a 27\;july\;1705 manuscript \emph{Puisque des personnes\ldots}
meant to be sent to Pierre Varignon (but sent to Jacques Lelong on
31\;july\;1705), Leibniz seeks to address one of the objections raised
by Michel Rolle concerning a case where l'H\^opital's rule needs to be
applied twice.  Here Leibniz writes:
\begin{quote}
Et ce cas des evanouissants ou naissants est si pres du cas dont il
s'agit qu'il n'en diff\`ere d'aucune \emph{grandeur assignable}, comme
il est manifeste dans la figure 130 du trait\'e de l'Analyse des
infinitesimales, ou [i.e., o\`u] la distance entre~$B$ et~$b$ est
moindre qu'aucune qu'on puisse \emph{assigner}.''  (Leibniz
\cite{Le05b}, 1705; emphasis added)
\end{quote}
The passage furnishes a succinct summary of Leibniz's take on
infinitesimals:
\begin{enumerate}
\item
the infinitesimal distance~$Bb$ is smaller than any assignable
magnitude;
\item
thus~$Bb$ is itself \emph{inassignable};
\item
the points~$B$ and~$b$ are specific points in an illustration from
l'Hospital's book;
\item
thus the distance~$Bb$ is (not a sequence but) a number.
\end{enumerate}
This apparently straightforward analysis of Leibniz's notion of
infinitesimal is resisted by a number of modern scholars.  Over three
centuries after Leibniz published his first article on infinitesimal
calculus in \emph{Acta Eruditorum} (\cite{Le84}, 1684), scholars are
still debating the nature of Leibnizian infinitesimals.  There are two
main interpretations of Leibnizian infinitesimals in the current
literature: one associated with historian Henk Bos and the other, with
philosopher Hid\'e Ishiguro.

\subsection{Bos on two approaches}
\label{s11}

Bos writes:
\begin{quote}
Leibniz considered two different approaches to the foundations of the
calculus; one connected with the classical methods of proof by
`exhaustion', the other in connection with a law of continuity.  (Bos
\cite{Bos}, 1974, p.\;55)
\end{quote}
Bos goes on to describe the second method as follows: ``The chief
source for Leibniz's second approach to the justification of the use
of `fictitious' infinitesimals in the calculus is a
manuscript\ldots{}'' \cite[p.\;56]{Bos}.  Here Bos is referring to
Leibniz's manuscript \emph{Cum Prodiisset} (\cite{Le01c}, 1701).

Leibniz's second approach mentioned by Bos is based on a law of
continuity, explained by examples as follows:
\begin{quote}
[I]n the case of intersecting lines, for instance, arguments involving
the intersection could be extended (by introducing an ``imaginary''
point of intersection and considering the angle between the lines
``infinitely small'') to the case of parallelism; also arguments about
ellipses could be extended to parabolas by introducing a focus
infinitely distant from the other, fixed, focus. \cite[p.\;57]{Bos}
\end{quote}
Bos grants that Leibnizian infinitesimals are fictional:
\begin{quote}
[Leibniz] had to treat the infinitesimals as `fictions' which need not
correspond to actually existing quantities, but which nevertheless
\emph{can be used} in the analysis of problems.
(\cite[pp.\;54--55]{Bos}; emphasis added)
\end{quote}

\subsection{Pure fictions and logical fictions}
\label{s13}

Such usable fictions could be termed \emph{pure fictions}; see
Sherry--Katz (\cite{14c}, 2014) in \emph{Studia Leibnitiana}.  Bos'
position is largely endorsed by Jesseph (\cite{Je15}, 2015).  The
question therefore is not \emph{whether} Leibnizian infinitesimals are
fictions, but rather \emph{what kind of} fictions: pure fictions or
logical fictions; see Section~\ref{s12}.

The matter of Appendix\;2 in (Bos \cite{Bos}, 1974) was dealt with in
Katz--Sherry (\cite{13f}, 2013, Section\;11.3, pp.\;606--608) and Bair
et al.\;(\cite{17b}, 2017, Section\;2.7, p.\;204); see also
Section~\ref{s37} here.

For a comparison of Leibniz's \emph{law of continuity} and the
\emph{transfer principle}%
\footnote{See note~\ref{f8} for a summary concerning transfer.}
of Abraham Robinson's framework for analysis with infinitesimals see
Katz--Sherry (\cite{12e}, 2012) as well as Section~\ref{s25}.  Bos'
position is largely endorsed by Paolo Mancosu; see Section~\ref{s24b}.
Marc Parmentier similarly sees two separate techniques in Leibniz's
work, including \emph{De Quadratura Arithmetica}:
\begin{quote}
Sa structure binaire se manifeste \'egalement dans l'oppo\-sition
entre deux types de m\'ethodes utilisables pour r\'ealiser des
quadratures.  Leibniz pr\'esente la premi\`ere comme un amendement,
d\^ument fond\'e et d\'emontr\'e, de la m\'ethode des indivisibles.%
\footnote{\label{f2}When Leibniz spoke of the method of indivisibles
in DQA, he sometimes had in mind a broader method of traditional
geometry originating with Archimedes, namely the technique of
exhaustion.  Thus, referring to his method \emph{not} exploiting
infinitesimals, Leibniz wrote: ``Quare methodo indivisibilium quae per
spatia gradiformia seu per summas ordinatarum procedit, ut severe
demonstrata licebit.'' (Parmentier's translation: ``Voil\`a ce qui
permettra de faire de la m\'ethode des indivisibles et de l'usage des
espaces gradiformes soit des sommes des ordonn\'ees qui en sont
l'apanage, une m\'ethode et un usage rigoureusement d\'emontr\'es''
\cite[p.\;63]{Le04b}.)  This method differs from the ``direct method''
which exploits infinitesimals.  See also note~\ref{f25}.}
Quant \`a l'autre m\'ethode, elle est fond\'ee sur les infiniment
petits.  L'enjeu de la \emph{Quadratura} est d\`es lors d'\'etablir
leur \'equivalence.  (Parmentier \cite{Pa01}, 2001, p.\;278)
\end{quote}
See Section~\ref{s312} and Section~\ref{s4} for more details on
\emph{De Quadratura Arithmetica}.

\subsection{Extensions and predicates}
\label{s15}

There is no need necessarily to rely on the idea of \emph{extension}
when formalizing Leibniz's procedures exploiting infinitesimals.
Indeed, in Edward Nelson's approach to infinitesimal analysis, such
entities are found in the ordinary real line.  Extensions
like~$\R\hookrightarrow\astr$ (see note~\ref{f12}) are necessary only
if one wishes to formalize infinitesimals in the context of
Zermelo--Fraenkel set theory based on the language possessing a single
relation, namely the membership relation~$\in$.  If, following Nelson,
one allows for a richer language including also a unary (i.e.,
one-place) predicate \textbf{standard} (together with axioms governing
its interaction with the Zermelo--Fraenkel axioms), then
infinitesimals (defined as a nonstandard numbers smaller in absolute
value than all positive standard ones) can be found within the
ordinary real line (see Section~\ref{s5} and note~\ref{s32} for
details).

The need for a two-tier number system to account for infinitesimal
calculus was felt by philosophers of the Marburg school, Hermann Cohen
(1842--1918) and Paul Natorp (1854--1924).  They exploited the pair
intensive/extensive to describe such a number system, with
infinitesimals being \emph{intensive} (Cohen's student Ernst Cassirer,
while nominally endorsing the view, in practice went on to analyze
other calculi, instead).  However, the anti-infinitesimal sentiment
fueled by Cantor, Russell, and others at the time was too powerful and
the necessary mathematical tools not yet available to enable a
convincing formalisation of such ideas; see Mormann--Katz\;\cite{13h}
for details.  An important link between the Marburg neo-Kantians and
Robinson's school is Abraham Fraenkel; see Kanovei et
al.\;(\cite{18i}, 2018) for details.  Felix Klein's take on
infinitesimals was more positive than is generally known; see Bair et
al.\;(\cite{18a}, 2017).

\subsection{Infinite sets \emph{vs} infinite numbers in Leibniz}
\label{s14b}

Leibniz held, following Galileo, that infinite aggregates,
collections, multitudes, or totalities (what we may refer to as
\emph{sets} today) lead to contradiction; see e.g., Knobloch
(\cite{Kn12}, 2012).  Thus, in a letter to Bernoulli dated
22\;august\;1698, Leibniz wrote:
\begin{quote}
To be sure, several years ago I have proved that the multitude of all
numbers implies a contradiction, if [it is] taken to be a single
totality.%
\footnote{In the original: ``Sane ante multos annos demonstravi,
numerum seu multitudinem omnium numerorum contradictionem implicare,
si ut unum totum sumatur.''}
(Leibniz \cite{Le98}, 1698, p.\;535)
\end{quote}
On the other hand, Leibniz routinely used fictional infinite (and
infinitesimal) \emph{numbers} in his work; see e.g.,
Section~\ref{s24}, Katz--Sherry (\cite{12e}, 2012), and Bl\aa sj\"o
(\cite{Bl17c}, 2017).

\subsection{Physics, matter, and space}
\label{f8c}

Leibniz occasionally uses the term \emph{syncategorematic} in
discussing physics, as described by De Risi:
\begin{quote}
[Leibniz] endorsed the (quite non-Aristotelian) view that bodies are
infinitely divided \emph{in actu}, but in a purely `syncategorematic'
way.  This latter notion aimed at expressing the idea that there is no
final element in the division of matter (i.e. no point), even though
there are more divisions of bodily parts than can possibly be
expressed by any finite number.  (De Risi \cite{De19}, 2019; emphasis
in the original).
\end{quote}
This use of the term ``syncategorematic infinity'' refers to matter or
space, closely related to indefinite divisibility.%
\footnote{See also main text at note~\ref{f25d}.}

\subsection
{Ishiguro, logical fictions, and alternating quantifiers}
\label{s12}

An alternative to Bos's interpretation was developed by Ishiguro in
(\cite{Is90}, 1990, Chapter\;5).%
\footnote{The interpretation in question was developed in the second,
1990 edition of Ishiguro's book, and is not yet found in its first
edition (Ishiguro \cite{Is72}, 1972).}
Ishiguro argues that a term that seems to express a Leibnizian
infinitesimal does not actually designate, denote, or refer; rather,
it is a \emph{logical fiction} in the sense of Russell; see e.g.,
(Russell \cite{Ru19}, 1919, p.\;45).  Such a reading of Leibnizian
infinitesimals contrasts with the Bos--Mancosu reading as presented in
Section~\ref{s13} in terms of \emph{pure fictions}.

The distinction between pure fiction and logical fiction
does not always receive sufficient attention from Leibniz scholars.
Thus, Ohad Nachtomy provides the following summary: 
\begin{quote}
Arthur argues that, due to a \emph{syncategorematic} interpretation of
the infinitely small, by 1676 Leibniz could use infinitesimals in
calculations and avoid the mystery -- and indeed the contradictions --
that their would-be existence would involve\ldots{} Thus a continuous
whole can be treated as if it consists in an infinity of
infinitesimals; but although by such means one can represent truths,
\emph{there are not such things in reality as infinite wholes or
infinitely small parts}.  (Nachtomy\;\cite{Na14}, 2014; emphasis
added)
\end{quote}
However, this passage amounts to a banal claim that in ``reality''
there is no referent for Leibnizian infinitesimals, an assertion
agreed to by Bos, Mancosu, and Parmentier.  As a summary of the
syncategorematic position that it purports to be, Nachtomy's passage
is a failure.

Note that Russell himself never applied his concept of \emph{logical
fiction} to Leibnizian infinitesimals.  According to Ishiguro,
\begin{quote}
[Leibniz] is treating [infinitely small lines] as convenient
theoretical fictions because using signs which looks [sic] as if they
stand for quantities sui generis is useful.''  (Ishiguro \cite{Is90},
1990, p.\;84)
\end{quote}
Such ``signs which look as if they stand for quantities'' turn out to
conceal universal and existential quantifiers as follows.  Ishiguro
contends that Leibniz's continuum is Archimedean, as when she
emphasizes ``the importance that Leibniz attached to his claim that,
strictly speaking, there are only finite numbers and magnitudes''
\cite[p.\;99]{Is90}.  She clarifies the nature of her non-designating
claim in the following terms: ``we can paraphrase the proposition with
a universal proposition with an embedded existential claim''
\cite[p.\;87]{Is90}.  In conclusion,
\begin{quote}
Fictions [such as Leibnizian infinitesimals] are \emph{not entities}
to which we refer.  \ldots{} They are correlates of ways of speaking
which can be reduced to talk about more standard kinds of entities.
(\cite[p.\;100]{Is90}; emphasis added)
\end{quote}
Such fictions, which are \emph{not entities to which we can refer} in
Ishiguro's view, are exemplified by Leibnizian infinitesimals:
\begin{quote}
We saw that Leibniz believed that his language of infinitesimals was
rigorous, although there is only a syncategorematic infinitesimal.
\cite[p.\;96]{Is90}
\end{quote}

\subsection
{Of rabbits and snakes}
\label{f4b}

Ishiguro's reading has been endorsed by a number of Leibniz scholars.
While not all of them subscribe to her \emph{not entity} view,
Ishiguro's idea that Leibnizian infinitesimals are \emph{not} to be
interpreted as non-Archimedean quantities has been endorsed by Arthur
(\cite{Ar13}, 2013, p.\;554), Goldenbaum (\cite{Go08}, 2008, p.\;76,
note\;59), Gray (\cite{Gr15}, 2015, p.\;10), Levey (\cite{Le15}, 2015,
p.\;184), Nachtomy (\cite{Na09}, 2009 and \cite{Na14}, 2014), and
elsewhere.  Rabouin (\cite{Ra15}, 2015, note\;25, p.\;362) describes a
Leibnizian infinitesimal as a ``{}`syncategorematic' entity'' citing
Ishiguro, but points out on the same page that ``this arbitrariness
[in the choice of~$\varepsilon$] does not amount, in modern terms, to
a universal quantification (at least in classical first order logic),
which would be meaningless to Leibniz.''

Knobloch interprets Leibnizian infinitesimals as \emph{variable
quantities}; see (\cite{Kn90}, 1990), (\cite{Kn94}, 1994),
(\cite{Kn02}, 2002), and (\cite{Kn08}, 2008); see also
Section~\ref{s23c}.  

Breger's view of Leibniz's infinitesimal as a \emph{process} led him
to a particularly colorful metaphor:
\begin{quote}
Whoever is interested in the provability rather than in the art of
finding should not stare at the infinitely small magnitude like a
rabbit at the snake; he should take a closer look at the process of
ever-decreasing divisions.  (Breger \cite{Br08}, 2008, p.\;188)
\end{quote}
Breger's claim that to understand Leibnizian infinitesimals one should
``look at the process of ever-decreasing divisions" may not be easy to
reconcile with Leibniz's claim that one does \emph{not} reach
infinitesimals by a process of ever-decreasing divisions; see e.g.,
Section~\ref{s24b}.  Breger goes on to offer a sharp criticism of Bos'
position in \cite[pp.\;196--197]{Br08}.  However, Breger's criticism
begs the question since it is predicated on the logical fiction
hypothesis.  Namely, Breger writes:
\begin{quote}
[A]ccording to Bos there were two strategies by which to justify
Leibniz's recourse to infinitesimals: epsilontics and the principle of
continuity (Bos, 1974, 55--57). This distinction appears artificial,
for the principle of continuity is \emph{of\;course} also founded on
epsilontics: two magnitudes are equal if their difference is smaller
than any magnitude that can possibly be expressed.  \ldots{} In either
case the processual nature is the decisive point; it is of no great
import whether the process is described by means of \emph{epsilontics}
or with reference to the principle of continuity.  (Breger
\cite{Br08}, 2008, p.\;197; emphasis added).
\end{quote}
Breger's rabbit-and-snake metaphor contains the germs of a remarkable
admission, namely that some historians tend to become paralyzed when
faced with \emph{bona fide} infinitesimals.%
\footnote{The ``epsilontic'' theme is continued in Breger's 2017
article: ``The fact that the infinitely small (and the incomparably
small) magnitudes derived their justification from epsilontics was
simply self-evident.  Leibniz did not consider it necessary to explain
this \emph{in any depth}'' (Breger \cite{Br17}, 2017, p.\;78; emphasis
added).  One appreciates Breger's admission that there is no
\emph{in-depth} source in Leibniz for the ``epsilontic'' reading.  See
further in note~\ref{f24}.}

\subsection{Syncategorematic talk}

Some of these authors claim to find support for their view in
Proposition\;6 from Leibniz's unpublished text \emph{De Quadratura
Arithmetica}.  Their claims have been challenged in a recent detailed
textual study by Bl\aa sj\"o (\cite{Bl17a}, 2017); see also Knobloch
(\cite{Kn17}, 2017), Bl\aa sj\"o (\cite{Bl17b}, 2017), and
Section~\ref{s4}.

Ishiguro's interpretation of Leibnizian infinitesimals goes under the
name \emph{syncategorematic}:%
\footnote{Peter Geach pointed out the inappropriateness of using the
term \emph{syncategorematic} in this context: ``{}`Categorematic' and
`syncategorematic'\ldots{} are words used to describe (uses of) words
in a language; an infinite multitude, say, can no more be
syncategorematic than it can be pronominal or adverbial.  To be sure,
the confusion is explicable\ldots'' (Geach \cite{Ge67}, 1967,
p.\;41).}
\begin{quote}
\ldots talk of infinitesimals is, as [Leibniz] says,
\emph{syncategorematic} and is actually about `quantities that one
takes\ldots{} as small as is necessary in order that the error should
be smaller than the given error.'  (Ishiguro \cite{Is90}, 1990,
p.\;90; emphasis added)
\end{quote}
Note that here Ishiguro applies the term not to physics, matter, or
space (see Section~\ref{f8c}), but rather to an individual magnitude.
Used this way, the term is the catchword for the idea that Leibnizian
infinitesimals are ``signs which look as if they stand for
quantities'' but in reality signify concealed quantifiers.%
\footnote{On this reading, they are ghosts of departed
\emph{quantifiers}; cf.\;Bair et al.\;(\cite{17a}, 2017).  Thus,
``[i]f first-order differentials have absorbed a logical quantifier,
second-order differentials have absorbed two logical quantifiers.''
(Breger \cite{Br08}, 2008, p.\;194).}
The passage Ishiguro is alluding to actually provides a piece of
evidence \emph{against} her syncategorematic thesis; see
Section~\ref{s31}.

\subsection{Of Leibniz and Weierstrass}
\label{s14}

Bowdlerized accounts of Leibniz's position portraying him as a
proto-Weierstrassian are ubiquitous in the literature.  Thus,
commenting on Leibniz's observation that ``it is unnecessary to make
mathematical analysis depend on or to make sure that there are lines
in nature which are infinitely small in a rigorous sense in contrast
to our ordinary lines, or as a result, that there are lines infinitely
greater than our ordinary ones, etc.''  (Leibniz \cite{Le89}, 1989,
pp.\;542--543), \emph{editor} Loemker feels compelled to declare:
\begin{quote} 
If Leibniz had more clearly combined his conception of the
infinitesimal as a quantity to be taken at will as less than any
assignable quantity whatever with his own analysis of series and his
functional conception of the law of continuity, he should have been
led to the critical concept of limits upon which the calculus was
\emph{at\;last} theoretically grounded in the nineteenth century by
Weierstrass and Cauchy.%
\footnote{On Cauchy see note~\ref{f6}.}
(Editor Loemker commenting in \cite[note\;2, p.\;546]{Le89}; emphasis
added)
\end{quote}
The presentist view of the history of analysis as inexorably
progressing toward, and culminating, \emph{at last}, in the
Weierstrassian \emph{Epsilontik} was analyzed by Bair et
al.\;(\cite{17a}, 2017).  See also Hacking (\cite{Ha14}, 2014) on the
distinction between a butterfly and a Latin model for the development
of a science.

\subsection
{Arthur's endorsement of syncategorematic reading} 
\label{s18}

In 2013 Richard Arthur endorses Ishiguro's reading in the following
terms:
\begin{quote}
I take the position here (following Ishiguro 1990) that the idea that
Leibniz was committed to infinitesimals as actually infinitely small
entities is a misreading: his mature interpretation of the calculus
was fully in accord with the Archimedean Axiom.  Leibniz's
interpretation is (to use the medieval term) \emph{syncategorematic}:
Infinitesimals are fictions in the sense that the terms designating
them can be treated as if they refer to entities incomparably smaller
than finite quantities, but really stand for variable finite
quantities that can be taken as small as desired.  (Arthur
\cite{Ar13}, 2013, p.\;554)
\end{quote}
Here Arthur seeks to apply the qualifier \emph{syncategorematic} to a
Leibnizian infinitesimal, rather than either multitude, matter, or
space (see Sections~\ref{s14b} and~\ref{f8c}).  The Ishiguro--Arthur
(IA) syncategorematic thesis concerning Leibniz's infinitesimals has
gained wide acceptance in the literature; see Section~\ref{s12}.  For
more details on (Arthur \cite{Ar13}) see Section~\ref{s38}.

In 2014 Arthur endorses an allegedly \emph{non-referring} nature of
Leibnizian infinitesimals:
\begin{quote}
Just as the infinite is not an actually existing whole made up of
finite parts, so infinitesimals are not existing parts which can be
composed into a finite whole.  Borrowing a term from the Scholastics,
Leibniz called the infinite and the infinitely small
\emph{syncategorematic terms}: like `it' or `some' in a meaningful
sentence, they do not in themselves \emph{refer} to determinate
things, but can be used perfectly meaningfully in a specified context.
(Arthur\;\cite{Ar14}, 2014, p.\;81; emphasis on `syncategorematic
terms' in the original; emphasis on `refer' added)
\end{quote}

In 2015 Arthur renews his endorsement in the following terms:
\begin{quote}
Ishiguro (1990)\ldots{} was one of the first to argue that Leibniz can
allow for the success of treating the infinite and infinitely small
\emph{as if} they are entities (under certain conditions), and that it
is this that allows him to claim that mathematical practice is not
affected by whether one takes them to be real or not.%
\footnote{\label{f6}Quantifiers are alluded to on the same page in the
following terms: ``the justification is in terms that, after
\emph{Cauchy}, we would now express in terms of~$\varepsilon$
and~$\delta$'' \cite[p.\;146]{Ar15} (emphasis added).  For an analysis
of the error of attributing prototypes of~$(\varepsilon,\delta)$
alternating quantifier definitions to Cauchy see Bascelli et
al.\;(\cite{18e}, 2018).}
(Arthur \cite{Ar15}, 2015, p.\;146, note\;16)
\end{quote}
Ishiguro attributes such \emph{nonentity} syncategorematic
\emph{as-if} infinitesimals to Leibniz without restricting it to any
specific period of Leibniz's career.  Arthur and Levey acknowledge the
presence of infinitesimal \emph{entities} at least in the early
Leibniz.  Accordingly, they have modified Ishiguro's position to a
syncategorematic interpretation starting as early as 1676 (see Arthur
\cite{Ar13}, 2013, p.\;554), recognizing that there is a historical
development of Leibniz's mathematical insights and ideas, including
the notion of infinitesimals.

In 2018, Arthur's syncategorematic infinitesimal goes \emph{actual}.
In a chapter entitled `Leibniz's syncategorematic actual infinite,'
Arthur writes:
\begin{quote}
[T]o say that a magnitude is \emph{actually} infinitely small in the
syncategorematic sense is to say that no matter how small a magnitude
one takes, there is a smaller, but there are no \emph{actual}
infinitesimals.'' (Arthur \cite{Ar18}, 2018, p.\;155; emphasis added)
\end{quote}
So is this magnitude actual or not\;actual?  Arthur's desire to
incorporate the qualifier ``actual'' in his title leads him to comical
incoherence in his discussion of magnitudes.  He goes on to offer yet
another endorsement of the \emph{logical fiction} hypothesis in the
following terms: ``In geometry one may calculate with expressions
apparently denoting such entities, on the understanding that they are
fictions, standing for variable magnitudes that can be made
arbitrarily small\ldots'' (ibid., pp.\;155--156).

A re-evaluation of Leibniz's contribution to analysis was developed in
2012 (Katz--Sherry \cite{12e}), in 2013 (Katz--Sherry \cite{13f}), in
2014 (Bascelli et al.\;\cite{14a}, Sherry--Katz \cite{14c}), in 2016
(Bair et al.\;\cite{16a}), in 2017 (Bair et al.\;\cite{17b},
B\l{}aszczyk et al.\;(\cite{17c}), and elsewhere.  In an apparent
reaction to this work, Arthur wrote in 2018:
\begin{quote}
Certain scholars of the calculus have denied that the interpretation
of infinitesimals as syncategorematic was Leibniz's mature view, and
have seen them as \emph{fictions in a different sense}.  I shall not
mainly be concerned with that line of disagreement here, reserving a
detailed critique of such views for another occasion.''
(Arthur\;\cite{Ar18}, 2018, p.\;156; emphasis added).
\end{quote}
Meanwhile Arthur's 2019 texts \cite{Ar19}, \cite{Ar19b} do not contain
the reserved critique.

An IA-style reading of Cauchy's infinitesimal as a logical fiction
has been challenged in Borovik--Katz (\cite{12b}, 2012), Bair et
al.\;(\cite{13a}, 2013), Bascelli et al.\;(\cite{14a}, 2014), Bair et
al.\;(\cite{17a}, 2017), Bascelli et al.\;(\cite{18e}, 2018), and
elsewhere.

\section{Mathematical fictions}

We argue that the IA position to the effect that ``Fictions [such as
Leibnizian infinitesimals] are not entities to which we refer\ldots{}
They are correlates of ways of speaking which can be reduced to talk
about more standard kinds of entities'' (see Section~\ref{s12})
involves equivocation on the meaning of the term \emph{fiction}.

\subsection{Entities, nonentities, and referents}

To the extent that we have symbolism for mathematical concepts, we
can, under suitable conditions, \emph{refer} to them; to the extent
that such concepts have no referents, we can also assert that they are
``not entities to which we refer'', provided we take note of the fact
that there is no difference here between Leibnizian infinitesimals on
the one hand and e.g., unending decimal strings,%
\footnote{These were already developed by Simon Stevin at the end of
the 16th century; see Katz--Katz (\cite{12c}, 2012).}
on the other.  In this sense Ishiguro's attempt at quantifier-assisted
transcription of infinitesimals amounts merely to an attempt at
long-winded paraphrase of one variety of nonentity by another.%
\footnote{Ishiguro uses the techniques of modern mathematics like
quantifiers to develop her reading of Leibniz, but overlooks the fact
that merely referring to real numbers as \emph{standard}, as she does,
does not make them any less lacking in referent than infinitesimals,
from the modern viewpoint.}

We argue that already in his first publication on the calculus in 1684
and especially starting in the 1690s, Leibniz exploited fictional
infinitesimals not reducible to a quantifier paraphrase, and even made
it clear that they violate the Archimedean property (see
Sections~\ref{s21} and \ref{s22b}).

\subsection{Smaller than any given quantity}

Leibniz repeatedly defined infinitesimals as being smaller than any
given quantity.  IA read this as ``a universal proposition with an
embedded existential claim'' (Ishiguro\;\cite{Is90}, 1990, p.\;87),
namely, as an assertion involving alternating quantifiers (see
Section~\ref{s12}).  However, a more straightforward reading is to
interpret \emph{given quantities} as being \emph{assignable} and
infinitesimals as \emph{inassignable}, in the terminology Leibniz used
both in \emph{Cum Prodiisset} in 1701 (see Section~\ref{s24}) and in
his manuscript \emph{Puisque des personnes\ldots} in 1705 (see
Section~\ref{s1}).%
\footnote{\label{f4}The distinction assignable \emph{vs} inassignable
goes back to the distinction quanta \emph{vs} non-quanta in work of
Nicholas of Cusa (1401--1464), which also inspired Galileo's
distinction between \emph{quanta} and \emph{non-quanta} according to
Knobloch (\cite{Kn99}, 1999, p.\;89 and \cite{Kn19}, 2019).  For
details on assignable \emph{vs} inassignable see Section~\ref{s24}.}
Sometimes Leibniz also uses the terminology of \emph{incomparables}
for infinitesimals (see Section~\ref{s22}).

Leibniz famously denied that infinite \emph{totalities} can be viewed
as \emph{wholes}, but such a rejection does not necessarily extend, at
least in the mathematical realm, to infinite and infinitesimal
\emph{magnitudes} and \emph{quantities}, as discussed in
Section~\ref{s14b}.

\subsection{Variable quantities from Varignon to Knobloch}
\label{s23c}

Knobloch reads Leibniz' infinitesimals as being variable
quantities: 
\begin{quote}
Eventually, Leibniz adhered to `smaller than any given quantity' or
infinitely small that is to a completely consistent fruitful
definition of infinitely small.  An infinitely small quantity is a
variable quantity and can be described in terms of the Weierstrassian
epsilon-delta language\ldots{} It must be a variable quantity that can
be described in the Weierstrassian~$\varepsilon$-$\delta$-language:
smaller than any given quantity.  (Knobloch, \cite{Kn19}, 2019,
pp.\;2,\;7)
\end{quote}
Knobloch's interpretation echoes a related stance expressed by
Varignon in 1700--01, when he attempted to defend and clarify Leibniz'
infinitesimal calculus in response to attacks by Rolle.  In this
section, we will compare Knobloch's view with Varignon's.  

Rolle subsequently engaged Joseph Saurin in an a polemic that lasted
several years.  We will analyze the Rolle--Saurin exchange in
Section~\ref{s26}.  For the purposes of this section, we quote the
following criticism expressed by Rolle:
\begin{quote}
On reconno\^\i t d'abord que les effets des m\'ethodes qu'on propose
dans la nouvelle Analyse, sont toujours les m\^emes quand on substitue
des quantit\'es finies \`a volont\'e au lieu des Infiniment petits
$dx$ \&~$dy$: ce qui prouve que le succ\`es, bon ou mauvais, n'est
point attach\'e \`a \emph{l'infinie petitesse} qu'on suppose dans le
Syst\^eme.  (Rolle \cite{Ro03b}, 1703, pp.\;324--325; emphasis added)
\end{quote}
Here Rolle speaks of arbitrary finite quantities in terms similar to
Knobloch's (with the exception of the ``Weierstrassian
$\epsilon$-$\delta$''), the difference being that Rolle is
\emph{contrasting} such quantities with infinitesimals, whereas
Knobloch seeks to \emph{identify} them.  Rolle argued that the new
\emph{Syst\^eme} is superfluous since whatever it can achieve can
already be achieved by what he saw as finitist algebraic techniques
developed by Fermat and others; see Section~\ref{s27}.  Clearly both
sides in the Rolle--Saurin exchange understood Leibnizian
infinitesimals to be characterized by a property Rolle refers to as
\emph{l'infinie petitesse}.  To pursue his Weierstrassian thesis,
Knobloch would be forced to postulate that Leibniz was misunderstood
by his contemporaries on account of his infinitesimals (Ishiguro faced
a similar dilemma; see Section~\ref{s311}).

This series of exchanges was part of a flurry of debate at the Paris
Academy of Sciences between 1700 and 1705 focused on the viability of
the new Leibnizian calculus; see Blay (\cite{Bl86}, 1986) and Mancosu
(\cite{Ma89}, 1989).

Rolle attacked the idea of infinitesimals in several ways, including
protesting that they are sometimes treated as nonzero quantities and
sometimes as absolute zeroes, and in particular that from the equation
$dx + x = x$, one must conclude~$dx = 0$; see \cite{Ma89}.

In his defense of Leibniz' infinitesimals against these arguments,
Varignon appealed to Newton's \emph{Principia Mathematica} as the
source of imagery and rigor, extensively quoting the Scholium to
Lemma\;XI in Book I.\, According to Varignon, Rolle had failed to
understand that infinitesimals are actually variable quantities, not
fixed ones.  They decrease continually until they reach zero, but are
``considered only in the moment of their evanescence''
\cite[p.\;231]{Ma89}.  

Varignon states (in his correspondence with Bernoulli) that
differentials consist ``in being infinitely small and infinitely
changing until zero, in being nothing but \emph{quantitates
evanescentes, evanescentia divisibilia}, they will always be smaller
than any arbitrary given quantity'' (Bernoulli \cite{Be88}, 1988,
p.\;357).

Varignon goes on to explain how this can always be expressed instead
with an exhaustion argument, ``in the way of the ancients":
\begin{quote}
Indeed, whatever difference can be assigned between two magnitudes
which differ only by a differential it will always be possible, on
account of the continual and indefinite variability of this infinitely
small differential, and as on the verge of being zero, to find a
differential less than the given difference.  Which shows, in the way
of the ancients, that notwithstanding their difference these two
quantities can be taken to be equal.  (In Bernoulli\;\cite{Be88},
1988, p.\;357)
\end{quote}
It is suggestive that Varignon had to look to Newton for an
appropriate account of infinitesimals as being changing quantities.

\subsection{Leibniz's response to Varignon}
\label{s24b}

In a letter dated 28 november 1701, Varignon asked Leibniz to make a
precise statement on what should be undestood by \emph{infinitesimal
quantity}.  Leibniz's letter \cite{Le02} dated 2\;february\;1702
contains a response to Varignon.  Mancosu summarizes Leibniz's
response in the following three points:
\begin{quote}
(a) There is no need to base mathematical analysis on metaphysical
assumptions.  (b)\;We can nonetheless admit infinitesimal quantities,
if not as real, [then] as well-founded fictitious entities, as one
does in algebra with square roots of negative numbers. Arguments for
this position depended on a form of the metaphysical principle of
continuity.  Or (c) one could organize the proofs so that the error
will be always less than any assigned error.  (Mancosu \cite{Ma96},
1996, p.\;172)
\end{quote}
What is conspicuously absent is the Newton--Varignon definition of
infinitesimal as a variable quantity (see Section~\ref{s23c}).
Mancosu notes further:
\begin{quote}
In his letter [Leibniz] merged \emph{two different foundational
approaches}.  The first was related to the classical methods of proof
by exhaustion; the second was based on a metaphysical principle of
continuity.  (ibid.; emphasis added)
\end{quote}
Thus, Mancosu follows Bos in seeing both A-track and B-track methods
(see Section~\ref{s41}) in Leibniz.  The B-track method (mentioned in
item\;(b) of in Mancosu's summary) relies on what is referred to by
Bos \cite[p.\;55]{Bos} and Katz--Sherry \cite{13f} as the \emph{law of
continuity}.

As Leibniz wrote in the 2\;february\;1702 letter to Varignon (GM, IV,
91--95), this is a heuristic law to the effect that the rules of the
finite are found to succeed in the infinite, and conversely the rules
of the infinite apply to the finite:
\begin{quote}
Yet one can say in general that though continuity is something ideal
and there is never anything in nature with perfectly uniform parts,
the real, in turn, never ceases to be governed perfectly by the ideal
and the abstract\ldots{} (Leibniz as translated in
\cite[p.\;544]{Le89})
\end{quote}
Having formulated the basic distinction between the real and the
ideal, Leibniz proceeds to formulate his heuristic law:
\begin{quote}
\ldots{} and that \emph{the rules of the finite are found to succeed
in the infinite}, \ldots{}%
\footnote{\label{f8b}The passage omitted at this point is analyzed in
Section~\ref{s25b}.}
\emph{And conversely the rules of the infinite apply to the finite},
as if there were infinitely small metaphysical beings, although we
have no need of them, and the division of matter never does proceed to
infinitely small particles.  (ibid.; emphasis added)
\end{quote}
We will analyze some reactions to Leibniz's heuristic law in
Section~\ref{s25}.

\subsection{\emph{Fixe et determin\'ee}}
\label{s25c}

The pair of qualifiers \emph{fixe et determin\'ee} occurs in several
letters in the Leibniz--Varignon exchange.  The first letter in the
series is a 28\;november\;1701 letter from Varignon to Leibniz,
complaining about Jean Galloys (Gallois) in the following terms:
\begin{quote}
M. l'Abb\'e Galloys\ldots{} repand ici que vous avez declar\'e
n'entendre par differentielle ou Infinement [sic] petit, qu'une
grandeur \`a la verit\'e tres petite, mais cependant toujours
\emph{fixe et determin\'ee}, telle qu'est la Terre par raport au
firmament, ou un grain de sable par raport \`a la Terre\ldots{}
(Varignon \cite{Va01}, 1701, p.\;89)
\end{quote}
Thus, Galloys (and Varignon following him) use the pair \emph{fixe et
determin\'ee} to refer to a specific ``small'' assignable magnitude.

In a 2\;february\;1702 response (quoted in Section~\ref{s41}), Leibniz
speaks of \emph{common incomparables} which are still ordinary
assignable numbers.  Leibniz describes the latter as not being fixed
and determined, by which he means that they need to be made
arbitrarily small in an exhaustion-type argument.

In a subsequent letter dated 22\;march\;1702 from Varignon to
Bernoulli (that reached Leibniz in april), Varignon mentions that he
showed Leibniz's 2\;february\;1702 letter to P. Gouye, and describes
the latter's ``choleric'' reaction to a perceived change in Leibniz's
stance (GM IV 97).  Leibniz responds on 14\;april\;1702 as follows:
\begin{quote}
Je reconnois d'avoir dit quelque chose de plus dans ma lettre, aussi
estoit-il necessaire, car il s'agissoit d'\'eclaircir le memoire, mais
je ne crois pas qu'il y ait de l'opposition.  Si ce Pere [Gouye] en
trouve et me la fait connoistre, je tacheray de la lever.  Au moins
n'y avoit il pas la moindre chose qui d\^ut faire juger que
j'entendois une quantit\'e tres petite \`a la verit\'e, mais tousjours
\emph{fixe et determin\'ee}.  (Leibniz \cite{Le02b}, 1702; emphasis
added)
\end{quote}
Here Leibniz uses the pair of qualifiers \emph{fixe et determin\'ee}
in the sense used by the opponents of the calculus (Galloys and
Gouye), namely to refer to a specific ``small'' assignable number.
Leibniz denies that his incomparable (as opposed to \emph{common}
incomparables mentioned in the 2\;february\;1702 letter) is such a
number.

The qualifier \emph{common} is similarly used in ``Tentamen de motuum
coelestium causis'' (Leibniz \cite{Le89b}, 1689) in reference to the
type of incomparables used to justify infinitesimals by means of
assignable numbers:
\begin{quote}
I have assumed in the demonstrations incomparably small quantities,
for example the difference between two \mbox{\emph{common}} quantities
which is incomparable with the quantities themselves.  Such matters as
these, if I am not mistaken, can be set forth most lucidly in what
follows.  And then if someone does not want to employ infinitely small
quantities, he can take them to be as small as he judges sufficient to
be incomparable, so that they produce an error of no importance and
even smaller than any given [error]. 
(Leibniz as translated by Jesseph
in \cite{Je08}, 2008, p.\;227; emphasis added)
\end{quote}
Leibniz proceeds to give practical examples of common incomparables:
\begin{quote}
Just as the Earth is taken for a point, or the diameter of the Earth
for a line infinitely small with respect to the heavens, so it can be
demonstrated that if the sides of an angle have a base incomparably
less than them, the comprehended angle will be incomparably less than
a rectilinear angle, and the difference between the sides will be
incomparable with the sides themselves; also, the difference between
the whole sine, the sine of the complement, and the secant will be
incomparable to these differences.  (ibid.)
\end{quote}

\subsection{The passage on atomism}
\label{s25b}

The passage we omitted in Section~\ref{s24b} (see note~\ref{f8b})
reads as follows: ``as if there were atoms, that is, elements of an
assignable size in nature, although there are none because matter is
actually divisible without limit'' (Leibniz as translated by Loemker).
Loemker's translation of this passage is imprecise.  We reproduce the
original:
\begin{quote}
comme s'il y avait des atomes (c'est \`a dire des elemens assignables
de la nature), quoyqu'il n'y en ait point la matiere estant
actuellement sousdivis\'ee sans fin.  (Leibniz \cite{Le02}, 1702,
p.\;93; original spelling retained)
\end{quote}
This is the unique mention of \emph{atomes} in Leibniz's letter.

Modern readers may well be puzzled by Leibniz's aside on atomism.
Granted Leibniz's consistent opposition to atomism, what need is there
to interrupt a discussion of the principles of infinitesimal calculus
by an aside concerning physical atomism?  Why does Leibniz feel a
need, specifically in an infinitesimal context, to distance himself
from atomism, a doctrine familiar to schoolchildren today?

A possible answer lies in the 17th century battles -- theological and
otherwise -- over the doctrines of hylomorphism, transubstantiation,
and eucharist.  Atomism was thought of by the catholic hierarchy at
the time as contrary to canon law as codified at the Council of Trent
in 1551 (session 13, canon 2) and therefore heretical.  The infinitely
small, via the language of indivisibles exploited by Cavalieri and
others, were thought of as closely related to atomism.  Leibniz's
aside therefore may have constituted a defensive move.  For details
see Fouke (\cite{Fo92}, 1992), Bair et al.\;(\cite{18d}, 2018).  

Amir Alexander (\cite{Al15}, 2015) offers a different account of the
opposition to indivisibles in the 17th century.  Namely, the jesuits
saw Euclidean mathematics as an organizing principle that helps man
make order out of chaos (and in particular defeat the reformers).
They saw indivisibles as introducing a dangerous discordant note in
the otherwise \mbox{(near-)perfect} harmony of Euclid, and therefore
opposed them as a subversive reform.  The jesuits viewed indivisibles
as actual errors introduced into the heart of pristine geometry.
Indivisibles made geometry paradoxical, unreliable, and chaotic, the
very opposite of what they believed it must be.  See also Sherry
(\cite{Sh18}, 2018) and Alexander (\cite{Al18}, 2018).  Inspite of
their differences, Alexander and Sherry agree on the following: (i)
indivisibles were controversial in the 17th century; (ii) the
opposition emanated from powerful religious circles; (iii) the
opposition was a major factor in the decline of the Italian school of
geometry.%
\footnote{The case of James Gregory is particularly instructive.
Gregory studied under the indivisibilist Stefano degli Angeli at
Padua.  Gregory left Padua in 1668 shortly before degli Angeli's
religious order of the jesuats was banned by papal brief in the same
year.  Gregory's books were subsequently supressed in Italy.  The same
year also marked an abrupt stop to degli Angeli's output on
indivisibles.  For details see Bascelli et al.\;(\cite{18f}, 2018).}

In a similar vein, Leibniz distanced himself from the idea of material
\emph{indivisibles} while discussing the fictional nature of
infinitesimals in a 20\;june\;1702 letter to Varignon:
\begin{quote}
Entre nous je crois que Mons.\;de Fontenelle, qui a l'esprit galant et
beau, en a voulu railler, lorsqu'il a dit qu'il vouloit faire des
elemens metaphysiques de nostre calcul.  Pour dire le vray, je ne suis
pas trop persuad\'e moy m\^eme, qu'il faut considerer nos infinis et
infiniment petits autrement que comme des choses ideales ou comme des
fictions bien fond\'ees.  \ldots{} Il est que les substances simples
(c'est \`a dire qui ne sont pas des estres par aggregation) sont
veritablement \emph{indivisibles}, mais elles sont immaterielles, et
ne sont que principes d'action.  (Leibniz \cite{Le02c}, 1702, p.\;110;
emphasis added)
\end{quote}
Leibniz expressed similar sentiments in a 1716 letter to Dagincourt;
see Section~\ref{s22}.

\subsection
{Leibniz's heuristic law: from Rolle to Robinson}
\label{s25}
Leibniz addressed a letter to Jean-Paul Bignon in july\;1705.  Here
Leibniz summarizes some objections voiced against infinitesimal
calculus by Rolle at the \emph{Academie} in the following terms:
\begin{quote}
[Ces objections] reviennent \`a dire en effect qu'en maniant ce
nouveau Calcul des infinitesimales, on ne doit point avoir la
libert\'e d'y joindre les \emph{axiomes et operations} de la Geometrie
et de l'Analyse ancienne; qu'on ne doit point substituer
\emph{aequalibus aequalia}, qu'on ne doit point dire que de deux
quantit\'es \'egales les quarr\'es sont egaux aussi, et choses
semblables; \ldots{} (Leibniz \cite{Le05}, 1705, p.\;838; emphasis in
the original; emphasis on ``axiomes et operations'' added)
\end{quote}
Leibniz speaks dismissively of such objections (going as far as
describing them as \emph{des chicanes} \cite[p.\;839]{Le05}).
Nonetheless, Rolle's objections are valid.  Why do the newly
introduced numbers obey the same axioms and operations as those
governing the \emph{ancienne} geometry and analysis?  Why does the
squaring operation extend as expected?  We will deal with these
objections in more detail in Section~\ref{s26}.

\subsection{Rolle on Descartes and Fermat}
\label{s27}

The thrust of Rolle's critique of infinitesimal calculus was that the
latter was both unnecessary and plagued by error.  Rolle felt that
infinitesimal calculus was unnecessary because the problems it solves
are solved more easily with what he claimed were ordinary algebraic
techniques already available.  Rolle was specifically referring to the
work of Fermat:
\begin{quote}
Cependant M. Descartes luy-m\^eme dans une autre Lettre \`a M. Hardy
explique \& perfectionne la Methode de M. de Fermat.  II designe la
\emph{difference} des abscisses par un segment de ligne dans la
figure, \& il la designe encore par la lettre~$e$ dans le calcul,
comme l'avoit d\'eja fait M.\;de Fermat luy-m\^eme.  Outre cela il
suppose une droite qui rencontre la courbe en deux points, \& qui doit
devenir Tangente lorsque la \emph{difference} ind\'etermin\'ee des
abscisses est prise pour un \emph{zero absolu}.  Il poursuit selon les
Regles ordinaires de la Geometrie \& de l'Algebre, \& selon les
id\'ees de l'Auteur dont il explique la Methode.  (Rolle \cite{Ro03},
1703, p.\;2)
\end{quote}
Fermat actually used a capital letter~$E$ (rather than the
lower-case~$e$).  The clause ``qui doit devenir Tangente lorsque la
\emph{difference} ind\'etermin\'ee des abscisses est prise pour un
\emph{zero absolu}'' is open to interpretation; see e.g.,
note~\ref{f25b}, (Katz et al.\;\cite{13e}, 2013), (Bair et
al.\;\cite{18d}, 2018).

\subsection{Rolle, Saurin, and \emph{chicanes}}
\label{s26}

The issue of handling the squaring operation in the extended domain,
mentioned in Section~\ref{s25}, is a subtler problem than might appear
at first sight.  Leibniz repeatedly emphasizes that he is working with
a generalized notion of equality ``up to'' a negligible term; see
e.g., Leibniz's sentences [1] and [2] quoted in Section~\ref{s21}.  It
is possible to interpret Leibniz's comment as asserting an equality
between, say, a pair of infinitely close \emph{infinite} numbers~$H$
and~$H+\epsilon$, where~$\epsilon$ is infinitesimal, e.g.,
$\epsilon=\frac{1}{H}$.  If so, computing the squares of the two
numbers we obtain a difference of
\[
(H+\epsilon)^2-H^2=H^2+2H\epsilon+\epsilon^2-H^2=
2H\epsilon+\epsilon^2=2+\epsilon^2.
\]
Thus, the difference between the squares in this case is an
appreciable (non-infinitesimal) amount~$2+\epsilon^2$, and one can
reasonably ask whether Leibniz would consider the squares still equal
and under what circumstances.  The generalized equality is used, for
example, in the proof of Leibniz's rule~$d(xy)=xdy+ydx$ which involves
dropping the negligible term~$dxdy$.  This procedure was pertinently
criticized by Rolle in \cite[p.\;327]{Ro03b}.%
\footnote{\label{f12}This issue can be routinely clarified in
Robinson's framework in terms of the standard part function, in the
context of the hyperreal extension~\mbox{$\R\hookrightarrow\astr$}.
The subring~$\hr\subseteq\astr$ consisting of the finite elements
of~$\astr$ admits a map~$\st$ to~$\R$, known as \emph{standard part}.
The map~$\st\colon \hr\to\R$ rounds off each finite hyperreal number
to its nearest real number.  This enables one, for instance, to define
the derivative of~$t=f(s)$ as~$f'(s)=\st\left(\frac{\Delta t}{\Delta
s}\right)$ (for infinitesimal~$\Delta s\ne0$).  For details see e.g.,
Keisler (\cite{Ke86}, 1986), Katz--Sherry (\cite{12e}, 2012).  See
also notes~\ref{f8} and \ref{f25b}.}

Leibniz's comment on squaring quoted in Section~\ref{s25} was prompted
by the following comment by Rolle:
\begin{quote}
On y voit un second d\'egagement, et une troisi\'eme substitution; on
y quarre les deux membres de la formule~$S=\frac{xdy}{dx}$\, Ce qui
n'a encore est\'e pratiqu\'e dans la G\'eometrie transcendante, ni
indiqu\'e par aucune regle dans cette G\'eometrie.  (Rolle
\cite{Ro03}, 1703, p.\;33)
\end{quote}
Leibniz's comment on extending ``axioms and operations'' was prompted
by Rolle's objection to what he felt was a dubious procedure that
consists in 
\begin{quote}
citer des regles ordinaires qui ayent quelque rapport aux operations;
supposer qu'elles sont particulieres \`a l'Analyse des Inf.\;petits.
\cite[p.\;33--34]{Ro03}
\end{quote}
Rolle's objections were quoted by Joseph Saurin in his response; see
(Saurin \cite{Sa05}, 1705, p.\;248).  It is instructive to examine
Saurin's reaction to Rolle's objections:
\begin{quote}
mais nous serions tomb\'es-l\`a, dans une extraction de racines, non
moins ino\"uie dans toute la Geometrie transcendante, que la
permission que nous nous sommes donn\'ee d'\'elever la formule au
quarr\'e; tant il nous \'etoit impossible de resoudre le cas propos\'e
par M. Rolle, sans faire des suppl\'emens \`a nos d\'efectueuses
M\'ethodes. \, (Saurin\;\cite{Sa05}, 1703, pp.\;253--254)
\end{quote}
Saurin's sarcasm is palpable, but what about a response to Rolle's
objection?  Alas, none is forthcoming.  Instead, sarcasm turns to
\emph{ad hominem}:
\begin{quote}
Apr\'es cette vaine \& puerile discussion, o\`u m'ont jett\'e les
difficultez de M.\;Rolle; je suis oblig\'e pour mon honneur de
d\'eclarer icy aux G\'eometres que je sens toute la honte qu'il y a
\`a s'arr\^eter \`a des objections de cette nature.  Si je le fais,
c'est parce qu'elles servent \`a faire conno\^\i tre de plus en plus
quel est l'esprit de l'Auteur que je refute, etc.  (ibid., p.\;254)
\end{quote}
Saurin proceeds next to Rolle's objection regarding the extension of
rules:
\begin{quote}
\emph{citer des regles ordinaires qui ayent quelque rapport aux
operations; supposer qu'elles sont particulieres \`a l'Ana\-lyse des
Inf. Pet.}  \ldots{} Tous cela paroles jett\'ee en l'air, \& qui ne
prouve autre chose, sinon que les manieres de l'Auteur sont to\^ujours
les m\^emes.  (ibid.)
\end{quote}
Saurin was unable to appreciate Rolle's objection but in fact, Rolle's
objection was more poignant than those formulated three decades later
by George Berkeley (see Katz--Sherry \cite{13f}, 2013).  We will
analyze Rolle's objection further in Section~\ref{s29}.

\subsection{\emph{La R\'eforme}}
\label{f13}

As in Leibniz's comments on atomism and indivisibles (see
Section~\ref{s25b}), religious tensions seem just below the surface in
the Rolle--Saurin exchange.  Rolle appears to feel free to exploit
phrases like \emph{selon la r\'eforme} when referring to Leibnizian
calculus, even though its practitioners never used the term to
describe the new techniques.  The term \emph{r\'eforme} occurs at
least ten times in (Rolle \cite{Ro03}).

Meanwhile, Saurin had converted to catholicism barely a decade
earlier, and would not necessarily have appreciated Rolle's choice of
terminology, containing an allusion to the Reformation and the
Counter-Reformation.  These were developments of a recent past at the
time.

In one of his responses, Saurin alludes to Rolle's terminology, and
asks rhetorically: where is the reform in all this?  (with ``reform"
italicized): ``Y a-t-il l\`a quelque \emph{d\'eguisement}, quelque
\emph{suppl\'ement}, quelque \emph{reforme}?''  (Saurin \cite{Sa06},
1706, p.\;12).  And again: ``D'abord on remarquera que cette solution
\`a laquelle M. Rolle s'est principalement attach\'e, est non \emph{un
nouveau suppl\'ement, une nouvelle reforme}, ainsi qu'il l'appelle,
mais une solution de surcro\^\i t; \ldots'' (ibid.)

Rolle's move of imputing ideologically impure motives to the
pro-infinitesimal opposition is not without modern imitators.%
\footnote{See Section~\ref{f20} at note~\ref{f21}.}

\subsection{Transfer}
\label{s29}

As noted in Section~\ref{s26}, Rolle was asking for a justification of
Leibniz's heuristic law allowing one to extend rules from the ordinary
domain to the extended domain of the infinitesimal calculus.  In
modern terminology such justification is provided by the
\emph{transfer principle}.  Robinson analyzed Leibniz's heuristic law
as follows:
\begin{quote}
Leibniz did say, in one of the passages quoted above, that what
succeeds for the finite numbers succeeds also for the infinite numbers
and vice versa, and this is remarkably close to our transfer%
\footnote{\label{f8}The \emph{transfer principle} is a type of theorem
that, depending on the context, asserts that rules, laws or procedures
valid for a certain number system, still apply (i.e., are
``transfered'') to an extended number system.  Thus, the familiar
extension~$\Q\hookrightarrow\R$ preserves the property of being an
ordered field.  To give a negative example, the
extension~$\R\hookrightarrow\R\cup\{\pm\infty\}$ of the real numbers
to the so-called \emph{extended reals} does not preserve the property
of being an ordered field.  The hyperreal
extension~\mbox{$\R\hookrightarrow\astr$} (see note~\ref{f12})
preserves \emph{all} first-order properties, e.g., the
formula~$\sin^2x + \cos^2 x =1$ which remains valid for all
hyperreal~$x$, including infinitesimal and infinite values
of~$x\in\astr$.}
of statements from~$\R$ to~$\astr$ and in the opposite direction.
(Robinson \cite{Ro66}, 1966, p.\;266)
\end{quote}
Rolle was unwilling to accept the validity of such transfer.  There is
little in the responses he received that could have satisfied him,
given his rejection of infinitesimal quantities.

Leibniz suggests that infinitesimals ought to be treated as fictitious
entities, as one does in algebra with square roots of negative
numbers.  Thus the difference is not merely a distinct approach to
infinitesimals from Varignon's, but a broader difference in their
stances on mathematical formalism.

According to Mancosu, ``by considering the infinitesimals as
well-founded fictions, [Leibniz] was introducing a gap between the
formal apparatus and the referents'' (Mancosu \cite{Ma96}, 1996,
p.\;173).

\section{Evidence: incomparables}

We argue that Leibnizian texts tend to support the interpretation of
Bos and Mancosu over that of IA.

As a general comment, note that someone seeking to contest the IA
interpretation of Leibnizian infinitesimals as logical \emph{fictions}
might be tempted to assert that Leibnizian infinitesimals are not
fictional but \emph{real}.  Such a formulation may unwittingly entail
ontological commitments as to the \emph{reality} of infinitesimals.
However, one can reject Ishiguro's interpretation and still retain the
fictionalist interpretation of Leibnizian infinitesimals.  

To borrow Moigno's and Connes' terminology, one might say that an
infinitesimal is \emph{chimerical}.%
\footnote{See Schubring (\cite{Sc05}, 2005, p.\;456) and Bascelli et
al.\;(\cite{18e}, 2018, Section\;4.1) on Moigno; see Connes
(\cite{Co04}, 2004, p.\;14) and Kanovei et al.\;(\cite{13c}, 2013,
Section\;8.2, p.\;287) on the views of Connes.}
It does not however follow that they are \emph{logical} chimeras in
the IA sense.

\subsection{Leibniz's syncategorematic passage}
\label{s31}

Leibniz wrote in 1702:
\begin{quote}
Cependant il ne faut point s'imaginer que la science de l'infini est
degrad\'ee par cette explication, \& reduite \`a des fictions; car il
reste toujours un infini syncategorematique, comme parle l'Ecole \& il
demeure vray par exemple, que 2 est autant que 1/1+1/2+1/4+1/8+
1/16+1/32+\&c. Ce qui est une \emph{series} infinie dans laquelle
toutes les fractions, dont les Numerateurs sont l'unit\'e, \& les
denominateurs de progression Geometrique double, sont comprises \`a la
fois; quoy qu'on n'y employe toujours que des nombres ordinaires,
(Leibniz \cite{Le02}, 1702; emphasis in the original)
\end{quote}
Having made his remark concerning what he refers to as
syncategorematic infinity that involves only ordinary numbers, Leibniz
goes on to conclude:
\begin{quote}
\& quoy qu'on n[']y fasse point entrer aucune fraction
\mbox{\emph{infiniment petite}}, ou dont le denominateur soit un
nombre infini.  (ibid.; emphasis added)
\end{quote}
The plain meaning of the passage is that an infinitely small fraction
(or a fraction whose denominator is an infinite number) is \emph{not}
involved in syncategorematic infinity.  Thus Leibniz takes what he
refers to as \emph{infinitely small fractions} to be \emph{bona fide}
infinitesimals of track\;B type, in the terminology of
Section~\ref{s41}.

\subsection{Euclid, Definition V.4, and incomparables}
\label{s21}

Leibniz repeatedly made it clear that his system of magnitudes
involves a violation of the Archi\-medean property, viz., Euclid's
\emph{Elements}, Definition V.4; see e.g., the passage in Leibniz
(\cite{Le95b}, 1695, p.\;322) as quoted by Bos (\cite{Bos}, 1974,
p.\;14).  This definition is a technical expression of Leibniz's
distinction between assignable and inassignable quantities; see
Sections~\ref{s1} and~\ref{s24}.

The violation of V.4 appears directly to contradict the IA claim that
Leibniz was working with an Archimedean continuum.  Leibniz frequently
writes that his infinitesimals are \emph{useful fictions}; but it is
best not to understand them as logical fictions but rather as
\emph{pure fictions}; see Section~\ref{s11}.

Let us consider in more detail Leibniz's comment in his article
\emph{Responsio ad nonnullas difficultates a Dn.\;Bernardo
Niewentiit\ldots} on Euclid V.5 (numbered V.4 in modern editions),
which is a version of the Archimedean axiom:%
\footnote{See De Risi \cite{De16}, 2016, Section\;II.3 for more
details on Euclid's Definition V.4.}
\begin{quote}
[1] Furthermore I think that not only those things are equal whose
difference is absolutely zero, but also those whose difference is
incomparably small.  [2] And although this [difference] need not
absolutely be called Nothing, neither is it a quantity comparable to
those whose difference it is.  [3] It is so when you add a point of a
line to another line or a line to a surface, then you do not increase
the quantity.  [4] The same is when you add to a line a certain line
that is incomparably smaller.  [5] Such a construction entails no
increase.  [6] Now I think, in accordance with Euclid Book V def.\;5,
that only those homogeneous quantities one of which, being multiplied
by a finite number, can exceed the other, are comparable.  [7] And
those that do not differ by such a quantity are equal, which was
accepted by Archimedes and his followers.%
\footnote{In the original Latin: ``[1] Caeterum aequalia esse puto,
non tantum quorum differentia est omnino nulla, sed et quorum
differentia est incomparabiliter parva; [2] et licet ea Nihil omnino
dici non debeat, non tamen est quantitas comparabilis cum ipsis,
quorum est differentia.  [3] Quemadmodum si lineae punctum alterius
lineae addas, vel superficiei lineam, quantitatem non auges.  [4] Idem
est, si lineam quidem lineae addas, sed incomparabiliter minorem.
[5]~Nec ulla constructione tale augmentum exhiberi potest.  [6]
Scilicet eas tantum homogeneas quantitates comparabiles esse, cum
Euclide lib.~5 defin.~5 censeo, quarum una numero, sed finito
multiplicata, alteram superare potest.  [7] Et quae tali quantitate
non differunt, aequalia esse statuo, quod etiam Archimedes sumsit,
aliique post ipsum omnes.''  (Leibniz \cite{Le95b}, 1695, p.~322;
numerals [1] through [7] added)}
(translated from Leibniz \cite{Le95b}, 1695, p.\;322; numerals [1]
through~[7] added)
\end{quote}
Here Leibniz employs the term \emph{line} in the sense of what we
would today call a \emph{segment}.  In clause [3], Leibniz exploits
the classical example with indivisibles (adding a point to a line
doesn't change its length) so as to motivate a similar phenomenon for
incomparables in clause~[4] (adding an incomparably small line to a
finite line does not increase its quantity), namely his
\emph{transcendental law of homogeneity} (Leibniz \cite{Le10b}, 1710)
summarized in clause\;[5].%
\footnote{In reference to this passage, Breger claims that ``the
unassignable magnitudes are fictitious, they cannot be determined by
any construction'' (Breger \cite{Br08}, 2008, p.\;196), but fails to
deal with Leibniz's very next sentence concerning Euclid\;V.5 (V.4 in
modern editions).}

In clause~[6], Leibniz refers to \emph{homogeneous} quantities
satisfying Euclid's definition~V.5, i.e., the Archimedean axiom.  In a
follow-up clause~[7], Leibniz goes on to refer to `those
[quantities],' say~$Q$ and~$Q'$, that `do not differ by such a
quantity,' namely they do \emph{not} differ by a \emph{homogeneous}
quantity of the type mentioned in clause [6] (that would
\emph{satisfy} Euclid~V.5 relative to~$Q$ or~$Q'$).  Rather,~$Q$
and~$Q'$ differ by a quantity \emph{not satisfying} Euclid~V.5, i.e.,
a quantity which violates V.5 relative to~$Q$ and~$Q'$.  Leibniz views
such quantities as equal in the sense of a generalized relation of
equality governed by his \emph{law of homogeneity}; see Katz--Sherry
(\cite{12e}, 2012).

Leibniz referred to differences as in clause [1] as \emph{incomparably
small}.  Thus Leibniz is clearly dealing with an \emph{incomparably
small} difference
$Q-Q'$ which \emph{violates} Euclid V.5 relative to~$Q$ or~$Q'$.

\subsection{Letter to l'Hospital refutes reading by Arthur}
\label{s22b}

Leibniz is even more explicit about the fact that his
\emph{incomparables} violate Euclid~V.5 in his letter to l'Hospital
dated 14/24\;june\;1695:
\begin{quote}
J'appelle \emph{grandeurs incomparables} dont l'une multipli\'ee par
quelque nombre fini que ce soit, ne s\c cauroit exceder l'autre, de la
m\^eme fa\c con qu'Euclide la pris dans sa cinquieme definition du
cinquieme livre.%
\footnote{Translation: ``I use the term \emph{incomparable magnitudes}
to refer to [magnitudes] of which one multiplied by any finite number
whatsoever, will be unable to exceed the other, in the same way
[adopted by] Euclid in the fifth definition of the fifth book [of the
\emph{Elements}]'' (V.4 in modern editions).}
(Leibniz \cite{Le95a}, 1695, p.\;288; original spelling preserved;
emphasis in the original)
\end{quote}
In formulas, what Leibniz is saying is that magnitude~$\varepsilon$ is
incomparable with a magnitude~$r$ when the following formula is
satisfied:~$(\forall{}n)[n\varepsilon<r]$ (for finite~$n$).  Thus
Leibniz makes it clear that his incomparable magnitudes violate the
Archimedean property (a.k.a.\;Euclid\;V.4) relative to~$r$.  We will
analyze Breger's discussion of this passage in Section~\ref{s34}.

Arthur's claim in \cite[p.\;562]{Ar13} based on the very passage from
\emph{Responsio} quoted in Section~\ref{s21} that allegedly ``Leibniz
was quite explicit about \emph{this Archimedean foundation} for his
differentials as `incomparables'\,'' (emphasis added) is therefore
surprising.  Arthur does not provide any explanation for his claim but
rather merely reproduces the passage we analyzed in Section~\ref{s21},
goes on to cite additional passages from Leibniz, and then gets into a
discussion of Leibniz's 1684 article and other texts.

Arthur thus fails to explain his inference of an allegedly Archimedean
nature of the Leibnizian continuum from this passage.  Therefore we
can only surmise the nature of Arthur's inference, apparently based on
the reference to Archimedes himself in the passage quoted in
Section~\ref{s21}.  However, the term \emph{Archimedean axiom} for
Euclid V.4 was not coined until the 1880s (see Stolz \cite{St83},
1883), about two centuries after Leibniz.  

Thus, Leibniz's mention of Archimedes in \cite{Le95b} could not refer
to what is known today as the Archimedean \emph{property} or
\emph{axiom}.  Rather, Leibniz mentions an ancient authority merely to
reassure the reader of the soundness of his methods, or alludes to the
method of exhaustion.  Arthur's cryptic claim concerning the passage
mentioning Archimedes (i.e., that it is indicative of an allegedly
Archimedean foundation for the Leibnizian differentials) is misleading
and anachronistic.

Leibniz's 1695 letter to l'Hospital (involving a violation of Euclid
Definition V.4 by Leibniz's incomparables) is absent from Arthur's
bibliography.  We will analyze Ishiguro's comments on the letter in
Section~\ref{s23b}.

\subsection{Breger on letter to l'Hospital}
\label{s34}

In Section~\ref{s22b} we presented our analysis of the allusion to the
violation of Euclid V.4 in Leibniz's 14/24\;june\;1695 letter to
l'Hospital.  Breger gives a similar interpretation of this passage:
\begin{quote}
In a letter to L'H\^opital of 1695, Leibniz gives an explicit
definition of incomparable magnitudes: two magnitudes are called
incomparable if the one cannot exceed the other by means of
multiplication with [sic] an arbitrary (finite) number, and he
expressly points to Definition 5 of the fifth book of Euclid quoted
above.  (Breger \cite{Br17}, 2017, p.\;73--74).  
\end{quote}
Here Breger acknowledges that Leibnizian incomparable magnitudes
violate Euclid\;V.5 (V.4 in modern editions), i.e., that they are
non-Archimedean relative to ordinary ones.  Breger's position is
especially significant.  This is because he otherwise pursues an
interpretation close to the logical fiction hypothesis (see e.g.,
Section~\ref{f4b}), leading Arthur to write: 
\begin{quote}
When Breger says that there is no actual infinite in Leibniz's
mathematics, he is primarily concerned to deny the reading of the
actual infinite in Leibnizian mathematics as categorematic (as in
\emph{non-Archimedean construals} of the continuum and
infinitesimals), and I have no quarrel with him\ldots{} about this.''
(Arthur \cite{Ar18}, 2018, p.\;157; emphasis added)
\end{quote}
Contrary to Arthur's claim (and to Breger's position expressed
elsewhere), Breger does not deny but rather endorses such a
``non-Archi\-me\-dean construal'' in \cite[pp.\;73--74]{Br17}.

We briefly consider the possibility that what Breger might have meant
here is an interpretation of Leibnizian incomparable magnitudes as
functions or sequences tending to zero.  However, such an
interpretation is untenable for the following reason.  If a violation
of V.4 is attributable to e.g., a sequence tending to zero, then it
becomes nearly impossible for a system of magnitudes to avoid being in
violation of Euclid V.4.  Namely, as soon as one incorporates
magnitudes corresponding to, say, the ordinary rational numbers, by
density one will be able to choose a sequence tending to zero, and
thus detect a violation of Euclid\;V.4 in this sense.  The only
systems \emph{not} violating Euclid\;V.4 would be discrete systems
like~$\N$.  It seems clear that this is not the meaning Euclid had in
mind when he formulated Definition 4 of his book V.  Since Leibniz
refers explicitly to Euclid it seems also clear that such a discrete
system of magnitudes is not what Leibniz had in mind, for otherwise he
would have likely mentioned such a significant departure from Euclid's
intention.

In other words, the ordinary system of magnitudes in Euclid is clearly
meant to obey V.4 for otherwise Euclid would not have stated V.4 as a
definition.  On the other hand, it is clear that such a system cannot
be as restrictive as a discrete system like~$\N$.  With regard to
Leibniz, it is particularly clear that his system of ordinary
(assignable) magnitudes necessarily incorporates arbitrarily small
ones, since Leibniz believed in indefinite divisibility of matter (see
e.g., Section~\ref{s25b}).  Such a system of ordinary magnitudes
cannot satisfy V.4 if Definition V.4 is interpreted in terms of
sequences.  Thus the interpretation of magnitudes in terms of
sequences is at tension with both Euclid's and Leibniz's intention.
See also end of Section~\ref{s23b}.

\subsection{Ishiguro on l'Hospital}
\label{s23b}

The 1695 letter \cite{Le95a} from Leibniz to l'Hospital responding to
Nieuwentijt's criticism is cited in Ishiguro's bibliography but
Ishiguro erroneously describes Leibniz's criticism of Nieuwentijt here
as criticism of\ldots{} Leibniz's ally, l'Hospital:
\begin{quote}
It is important to realise however that in this letter Leibniz is
using de l'Hospital's own criterion to refute him.  De l'Hospital had
asserted both that higher differentials are not magnitudes and, if
after being multiplied by an infinite number the assumed quantity does
not become an ordinary magnitude, then it is not a magnitude at all.
It is a nothing.  Leibniz responded that if that is what de l'Hospital
believes, then he cannot at the same time claim that~$ddx$ and~$dxdx$
are not magnitudes, since they would, if multiplied by an infinite
number (``\emph{per numerum infinitum sed altiorem seu infinites
infinitum}") become ordinary magnitudes.  This is, however, \emph{not}
Leibniz's \emph{own} criterion, as he does not believe that there is
such a thing as an infinite number.%
\footnote{Ishiguro provides no sources to justify her claim that
``Leibniz does not believe that there is such a thing as an infinite
number.''  Nor does she pay attention to the distinction between
multitude and number; see Section~\ref{s14b}.}
He is on the contrary trying to explain what differentials and
quadratures are by spelling out the thought that leads to them in
terms of finite quantities, finite numbers, and Leibniz's concept of
`infinitely many' and of `incomparable.'  (He points out for example
that de l'Hospital is wrong to think that if~$dy$ is equal,~$dx$ would
also be equal.).  (Ishiguro \cite{Is90}, 1990, p.\;89)
\end{quote}
Ishiguro repeatedly attributes to l'Hospital what Leibniz describes as
Nieuwentijt's errors.  Nieuwentijt is not mentioned at all on
Ishiguro's page\;89.

What are we to make of Ishiguro's command of the details of the
alignment of forces among Leibniz's contemporaries with regard to
infinitesimal calculus?  The shoddiness of her command of such details
undermines the credibility of her sweeping claims to the effect that
Leibniz was allegedly misunderstood by his contemporaries like
Bernoulli and l'Hospital (see Section~\ref{s311}), particularly in
view of the fact that Leibniz specifically endorses l'Hospital's
approach; see Section~\ref{s1}.

No scholar of ancient Greece has yet stepped forward to give a
syncategorematic reading of Euclid's Definition V.4.  It seems
reasonable to assume that Leibniz's understanding of Euclid's
Definition\;V.4 and its negation was similar to that of modern
scholars; see e.g., De Risi (\cite{De16}, 2016).  Therefore, to
account for Leibniz's 1695 texts analyzed in Sections~\ref{s21} and
\ref{s22b}, advocates of the logical fiction approach would have to
extend Ishiguro's hypothesis that Leibniz was misunderstood by his
contemporary scholars to apply to modern scholars of Euclid, as well.

\subsection
{\emph{Nouveaux Essais sur l'Entendement Humain}}
\label{s23}

In his 2014 book, Arthur makes the following claim:
\begin{quote}
Having reached this conclusion in 1676, [Leibniz] holds it from then
on: `there is no infinite number, nor infinite line or other infinite
quantity, if these are taken to be genuine wholes.' (NE 157) There is
an actual infinite, but it must be understood
\emph{syncategorematically}, etc.  (Arthur \cite{Ar14}, 2014, p.\;88)
\end{quote}
The same syncategorematic claim, based on the same Leibnizian passage,
is reproduced four years later in (Arthur \cite{Ar18}, 2018, p.\;161).
However, a careful examination of the evidence leads one to the
opposite conclusion from Arthur's.  Arthur's reference (NE 157) is an
English translation of Leibniz's treatise \emph{Nouveaux Essais sur
l'Entendement Humain}.  Here a fictional character named Th\'eophile
argues as follows:
\begin{quote}
\textbf{Th\'eophile:} A proprement parler, il est vrai qu'il y a une
infinit\'e de choses, c'est-\`a-dire qu'il y en a toujours plus qu'on
puisse assigner.  Mais il n'y a point de nombre infini ni de ligne ou
autre quantit\'e infinie, si on les prend pour des v\'eritables touts,
comme il est ais\'e de d\'emontrer.  Les \'ecoles ont voulu ou d\^u
dire cela, en admettant un infini syncat\'egor\'ematique, comme elles
parlent, et non pas l'infini cat\'egor\'ematique.  (Leibniz
\cite{Le04}, 1704, p.\;113)
\end{quote}
This passage occurs in Chapter 17 (of Book II) entitled \emph{De
l'Infinit\'e.}  Arthur reproduces the first two sentences of the above
passage, but fails to report the outcome of the discussion between
Th\'eophile and another fictional character, Philal\`ethe.  The above
preliminary comment by Th\'eophile is in response to the opening
comment by Philal\`ethe:
\begin{quote}
\textbf{Philal\`ethe: 1} Une notion des plus importantes est celle du
\emph{fini} et de l'\emph{infini}, qui sont regard\'ees comme des
modes de la quantit\'e.  (ibid.)
\end{quote}
A disagreement soon emerges between the two characters.  While
Philal\`ethe views the finite and infinite as ``des modifications de
l'\'etendue et de la dur\'ee,'' Th\'eophile insists that ``la
consid\'eration du fini et infini a lieu partout o\`u il y a de la
grandeur et de la multitude.''  On the latter view, the infinite is an
attribute of magnitude and multitude.  It is not \emph{not} an
attribute of extension (i.e., continuum or space) and time, as
Philal\`ethe argues.  When Philal\`ethe again attempts to connect the
infinite to \emph{space}, Th\'eophile provides a detailed rebuttal,%
\footnote{\label{f25d}This connects with Leibniz's syncategorematic
views of physical space, as mentioned by de Risi; see
Section~\ref{f8c}.}
and concludes as follows:
\begin{quote}
Mais on se trompe en voulant s'imaginer un espace absolu qui soit un
tout infini compos\'e de parties, il n'y a rien de tel, c'est une
notion qui implique contradiction, et ces touts infinis, et leurs
oppos\'es infiniment petits, ne sont de mise que dans le \emph{calcul
des g\'eom\`etres}, tout comme les racines imaginaires de l'alg\`ebre.
(Leibniz \cite{Le04}, 1704, p.\;114; emphasis added)
\end{quote}
Thus according to Th\'eophile, in space, \emph{bona fide}
infinitesimals are impossible; but in the \emph{calculations of
geometers}, they do have a place.  Th\'eophile's suggestion that
infinitesimals are possible in calculation on par with imaginary
numbers is cogent; see also Section~\ref{s45}.  Arthur fails to
reproduce this crucial passage, which constitutes a piece of evidence
against the IA \emph{logical fiction} hypothesis, since there does not
exist a syncategorematic paraphrase of imaginary numbers as such
logical fictions.  Thus the very chapter 17 of \emph{Nouveaux Essais}
from which Arthur quotes in support of the IA thesis actually
furnishes evidence \emph{against} it.

\subsection{Arnauld--Leibniz exchange}

In (\cite{Ar15}, 2015), Arthur equivocates on the exact meaning of the
\emph{syncategorematic} claim with regard to its implications for
Leibniz's infinities and infinitesimals, but whether or not his
reading is identical to Ishiguro's, Arthur does not seek to
differentiate his reading from Ishiguro's and on the contrary
repeatedly endorses Ishiguro's \emph{logical fiction} reading; see
Section~\ref{s18}.  No less an authority than Knobloch read the text
(Arthur \cite{Ar15}, 2015) as creating an impression that the
syncategorematic reading of Leibniz's mathematical infinitesimals
finds support in the Arnauld--Leibniz exchange in 1687.  Thus, in his
review of Arthur's text \cite{Ar15} for MathSciNet, Knobloch notes
that
\begin{quote}
[Richard] Arthur essentially bases his [syncategorematic] de\-ductions
on Leibniz's correspondence with [Antoine] Arnauld and [Burchard] de
Volder.  (Knobloch \cite{Kn15}, 2015)
\end{quote}
What exactly is the basis, allegedly deriving from such
correspondence, for the IA thesis?  Arnauld being more influential
than de Volder, we will focus on the Arnauld--Leibniz exchange.
Arthur cites the following passage from Leibniz's letter to Arnauld:
\begin{quote}
I believe that where there are only beings by aggregation, there will
not in fact be any real beings; for any being by aggregation
presupposes beings endowed with a true unity, because it derives its
reality only from that of its constituents. It will therefore have no
reality at all if each constituent being is still a being by
aggregation, for whose reality we have to find some further basis,
which in the same way, if we have to go on searching for it, we will
never find.%
\footnote{In the original: ``[J]e croy que l\`a, o\`u il n'y a que des
estres par aggregation, il n'y aura pas m\^eme des estres reels; car
tout estre par aggregation suppose des estres dou\'es d'une veritable
unit\'e, parcequ'il ne tient sa realit\'e que de celle de ceux dont il
est compos\'e, de sorte qu'il n'en aura point du tout, si chaque estre
dont il est compos\'e est encor un estre par aggregation, ou il faut
encor chercher un autre fondement de sa realit\'e, qui de cette
maniere s'il faut tousjours continuer de chercher ne se peut trouver
jamais'' (Leibniz \cite{LA}, pp.\;91--92).}
(Leibniz to Arnauld, 30\;april\;1687 as translated in
\cite[p.\;152]{Ar15})
\end{quote}
The reader may well wonder what, if anything, this has to do with
mathematical infinitesimals.  Indeed, the context of exchange between
Arnauld and Leibniz was the latter's views as detailed in his
``Discourse on Metaphysics'' dating from 1686.  Arnauld addressed a
letter to Leibniz on 4\;march\;1687, and Leibniz replied on
30\;april\;1687.  Arnauld and Leibniz are discussing the
\emph{metaphysics} related to the structure of \emph{matter}, rather
than anything related to mathematical infinitesimals, as is evident
from Leibniz's very next sentence:
\begin{quote}
J'accorde, Monsieur, que dans toute la nature corporelle il n'y a que
des machines (qui souvent sont anim\'ees) mais je n'accorde pas qu'il
n'y ait que des aggreg\'es de substances, et s'il y a des aggreg\'es
des substances, il faut bien qu'il y ait aussi des veritables
substances dont tous les aggreg\'es resultent.'  \cite[p.\;92]{LA}
\end{quote}
Thus, the Arnauld--Leibniz exchange is not concerned with the nature
of mathematical infinitesimals, contrary to Arthur's claim as reported
by Knobloch.  Leibniz's position on mathematical infinitesimals is
well known: it is not necessary to make mathematical analysis
dependent upon \emph{metaphysical} controversies
(cf.~Section~\ref{s14}):
\begin{quote}
``my intention was to point out that it is unnecessary to make
mathematical analysis depend on \emph{metaphysical} controversies or
to make sure that there are lines in \emph{nature} which are
infinitely small in a rigorous sense in comparison with our ordinary
lines\ldots{}'' (Leibniz as translated in \cite[p.\;86]{Is90};
emphasis added)
\end{quote}
Seeking an explanation of the nature of infinitesimals in Leibniz's
comments on the substantial status of being of aggregates of material
things, as Arthur attempted to do, is sheer obfuscation.

\subsection{Leibniz and Smooth Infinitesimal Analysis}
\label{s38}

Peckhaus mentions that Arthur compares Leibnizian infinitesimals with
nilpotent infinitesimals of Smooth Infinitesimal Analysis (SIA; see
e.g., Bell \cite{Be08}):
\begin{quote}
In Section\;3 Leibniz's conception is compared with\ldots{} Bell's
Smooth Infinitesimal Analysis (SIA) \ldots{} which has many points in
common with the Leibnizian approach.  (Peckhaus \cite{Pe13}, 2013).
\end{quote}
In fact the original title of Arthur's text was ``Leibniz's
syncategorematic infinitesimals, Smooth Infinitesimal Analysis, and
Newton's Pro\-position 6.''%
\footnote{A preprint of Arthur's article with this longer title can be
viewed at \url{http://u.math.biu.ac.il/~katzmik/arthur08.pdf}}
Peckhaus notes moreover that Arthur finds many points in common as
well as dissimilarities:
\begin{quote}
The author comes to the conclusion that despite many points in common,
Leibniz's syncategorematic approach to the infinitesimals and Smooth
Infinitesimal Analysis ``are by no means equivalent''.  (ibid.)
\end{quote}
Note that SIA depends crucially on a category-theoretic foundational
framework and on intuitionistic logic to enable nilpotency.  While we
welcome such a display of foundational pluralism on Arthur's part, we
also agree with both John L. Bell and Arthur that SIA is closer to
Nieuwentijt's approach to calculus than to Leibniz's.  

Note that the nilpotent infinitesimals of SIA are smaller in absolute
value than~$\frac{1}{n}$, signifying non-Archimedean behavior.  The
following question therefore arises.  Why does Arthur's voluminous
output on Leibnizian calculus systematically eschew readings related
to Robinson's framework and Nelson's foundational approach (see
Section~\ref{s15}), based as they are on classical set-theoretic
foundations and classical logic, and enabling straightforward
transcription of both Leibniz's assignable/inassignable dichotomy and
his infinitely many orders of infinitesimal and infinite numbers?  The
question was essentially posed five years ago in the \emph{Erkenntnis}
article \cite{13f} and still awaits clarification.  Meanwhile, Arthur
appears to avoid Robinsonian infinitesimals as zealously as Leibniz
avoided atoms and material indivisibles; see Section~\ref{s25b}.

\subsection{Leibniz, Hilbert, and Formalism}
\label{s22}

The driving force behind the IA interpretation of Leibnizian
infinitesimals seems to be a desire to deny them an ontological
reality that would smack of Platonism or mathematical realism.
Leibniz himself was clear on the matter:
\begin{quote}
Quand [nos ami] disput{\`e}rent en France avec l'Abb{\'e} Gallois
[i.e., Galloys], le P{\`e}re Gouge [i.e., Gouye]~\&{} d'autres, je
leur t{\'e}moignai, que je ne croyois point qu'il y e{\^u}t des
grandeurs v{\'e}ritablement infinies ni v{\'e}ritablement
infinit{\'e}simales, que ce n'{\'e}toient que des fictions, mais de
fictions utiles pour abr{\'e}ger~\&{} pour parler universellement,
comme les racines imiginaires dans l'Alg{\`e}bre, telles
que~$\sqrt[2]{(-1)}$ \ldots{} (Leibniz \cite{Le16}, 1716)
\end{quote}
Without choosing sides in the realist/antirealist debate, we note that
Leibnizian infinitesimals can be understood as fictions in a sense
close to the school of Formalism as developed in the 20th century
mathematics by David Hilbert (see e.g., Hilbert \cite{Hi26}, 1926,
p.\;165) and others.

Thus, Abraham Robinson as a formalist distanced himself from platonist
and foundationalist views in the following terms:
\begin{quote}
[M]athematical theories which, allegedly, deal with infinite
totalities do not have any detailed \ldots{} reference.  (Robinson
\cite{Ro75}, 1975, p.\;42)
\end{quote}
Robinson explicitly linked the approaches of Leibniz and Hilbert in
the following terms: ``Leibniz's approach is akin to Hilbert's
original formalism, for Leibniz, like Hilbert, regarded infinitary
entities as ideal, or fictitious, additions to concrete Mathematics''
(Robinson \cite{Ro67}, 1967, pp.\;39--40).

\subsection{A Robinson--G\"odel exchange}
\label{s37}

Bos comments as follows on the connection between the frameworks of
Leibniz and Robinson:
\begin{quote}
The creation o<f non-standard analysis makes it necessary, according to
Robinson, to supplement and redraw the historical picture of the
development of analysis \ldots{} it is understandable that for
mathematicians who \emph{believe that [the] present-day standards [of
mathematical rigor] are final}, nonstandard analysis answers
positively the question whether, after all, Leibniz was right.
However, I do not think that being ``right" in this sense is an
important aspect of the appraisal of mathematical theories of the
past.  (Bos \cite{Bos}, 1974, p.\;82; emphasis added)
\end{quote}
This criticism of Robinson by Bos is predicated on the assumption that
Robinson believed that the present-day standards of mathematical rigor
are ``final.''  However, the attribution to Robinson of such naive
realist views concerning the finality of this or that piece of
mathematics is unsourced and unjustified, as we argue in this section
based on the Robinson--G\"odel correspondence.

In a 23\;august\;1973 letter to Kurt G\"odel, Robinson refers to his
posthumously published paper (\cite{Ro75}, 1975, presented at Bristol
in 1973) on progress in the philosophy of mathematics.  In this paper,
Robinson expresses formalist views.  The paper was enclosed with the
23\;august\;1973 letter to G\"odel.  Robinson's letter is a follow-up
on discussions that took place between him and G\"odel during
Robinson's visit to the Institute for Advanced Study (15--18
august\;1973).

In the letter Robinson writes: ``I am distressed to think that you
consider my emphasis on the model theoretic aspect of Non Standard
Analysis wrongheaded'' and goes on to describe himself as a ``good
formalist''.  Robinson then expresses the sentiment that G\"odel is
``bound to disagree" with the paper, possibly due to G\"odel's realist
views; see G\"odel (\cite{Go}, 2003; particularly the introduction by
M.~Machover) and Dauben (\cite{Da95}, 1995, pp.\;268--269).

Without getting into a discussion of the nature and extent of Goedel's
realist views, what we wish to highlight is the inaccuracy of
attributing such views to Robinson.  In his formalist views Robinson
was close to both Leibniz and Hilbert.

\section{Evidence: inassignables}

\subsection{Inassignable~$dx$ and assignable~$(d)x$}
\label{s24}

In addition to ratios of inassignable differentials such
as~$\frac{dy}{dx}$, Leibniz also considered ratios of ordinary values
which he denoted $(d)y$ and~$(d)x$, so that~$\frac{(d)y}{(d)x}$ would
be what we call today the derivative.  Here~$dx$ and~$(d)x$ are
distinct entities since Leibniz describes them as respectively
\emph{inassignable} and \emph{assignable} in \emph{Cum Prodiisset}
\cite{Le01c}:
\begin{quote}
[A]lthough we may be content with the assignable quantities
$(d)y$,~$(d)v$,~$(d)z$,~$(d)x$, etc., \ldots{} yet it is plain from
what I have said that, at least in our minds, the unassignables
[\emph{inassignables} in the original Latin]~$dx$ and~$dy$ may be
substituted for them by a method of supposition even in the case when
they are evanescent; \ldots{} (Leibniz as translated in Child
\cite{Ch}, 1920, p.\;153)
\end{quote}
Leibniz used similar terminology in his manuscript \emph{Puisque des
personnes}\ldots (\cite{Le05b}, 1705); see Section~\ref{s1}.  In
Leibniz (\cite{Le95b}, 1695), one similarly finds:
\begin{quote}
\ldots Nous voyons par l\`a que nous pouvons faire comme si le calcul
diff\'erentiel ne concernait que des quantit\'es ordinaires, m\^eme
s'il faut en rechercher l'origine dans les \emph{inassignables} pour
rendre compte des termes qui sont \'elimin\'es ou se d\'etruisent.
(Leibniz as translated by Parmentier in \cite{Le89}, 1989, p.\;336;
emphasis added).
\end{quote}
Meanwhile, on the IA reading,~$dx$ and~$(d)x$ should be identical,
both denoting ordinary assignable values (perhaps equipped with a
hidden quantifier or placed in a sequence).  The distinction between
differentials~$dx$ and~$(d)x$, extensively commented upon by Bos
(\cite{Bos}, 1974), is an indication that Leibniz exploits
differentials as pure fictions.

This is particularly significant since in \emph{Cum Prodiisset}
Leibniz is actually \emph{doing} calculus (thus, he proves the product
law for differentiation -- Leibniz's rule -- relying on the
transcendental law of homogeneity; see Katz--Sherry \cite{12e}, 2012),
rather than merely \emph{speculating} about it.  Breger wrote:
\begin{quote}
It has often been noted that Leibniz's \emph{verbal} descriptions of
infinitesimal magnitudes vary or even appear \emph{incoherent}
\ldots{} But in his \emph{use} of them Leibniz is in fact being quite
\emph{clear and explicit}; his view of infinitesimals appears not to
have altered since the beginning of his Hannover period or a few years
later. It is not sufficient to study Leibniz's \emph{verbal}
descriptions of infinitesimal magnitudes in isolation; they need to be
interpreted in connection with their \emph{mathematical usage}.
(Breger\;\cite{Br08}, 2008, p.\;185) (emphases added)
\end{quote}
We disagree with Breger's claim of alleged incoherence of Leibniz's
verbal descriptions, but we agree concerning the need to focus on
mathematical usage.

\subsection{Characteristic triangle}
\label{s25d}

\begin{figure}[ht]
\begin{center} 
\begin{tikzpicture} %
    %
    %
    %
    \node[circle,draw,fill,inner sep=0em,minimum size=1.0mm,black,opacity=1.0,label=above right:$T$] (T) at (0pt,75pt){};
    \node[circle,draw,fill,inner sep=0em,minimum size=1.0mm,black,opacity=1.0,label=above right:$A$] (A) at (0pt,-40pt){};
    \node[circle,draw,fill,inner sep=0em,minimum size=1.0mm,black,opacity=1.0,label=above right:$B$] (B) at (0pt,-130pt){};
    \node[circle,draw,fill,inner sep=0em,minimum size=1.0mm,black,opacity=1.0,label=above right:$(B)$] (Bparen) at (0pt,-180pt){};
    \draw[line width=0.75pt] (T) -- ($(Bparen)+(0pt,-10pt)$);
    %
    %
    \node[circle,draw,fill,inner sep=0em,minimum size=1.0mm,black,opacity=1.0,label=above right:$C$] (C) at (113pt,-130pt){};
    \node[circle,draw,fill,inner sep=0em,minimum size=1.0mm,black,opacity=1.0,label=above right:$D$] (D) at (113pt,-180pt){};
    \draw[line width=0.75pt,dashed] (C)--(D);
    %
    %
    \node[circle,draw,fill,inner sep=0em,minimum size=1.0mm,black,opacity=1.0,label=above right:$(C)$] (Cparen) at (140pt,-180pt){};
    \draw (A) parabola ($(Cparen)+(4.5pt,-8.5pt)$);
    %
    %
    %
    \draw (B)--(C);    
    \draw (Bparen)--(Cparen);
    %
    %
    \draw [line width=0.25pt,opacity=0.5] (T)--(Cparen);
  \end{tikzpicture}  
 
\end{center}
\caption{Leibniz's tangent line~$TC(C)$}
\label{f321}
\end{figure}
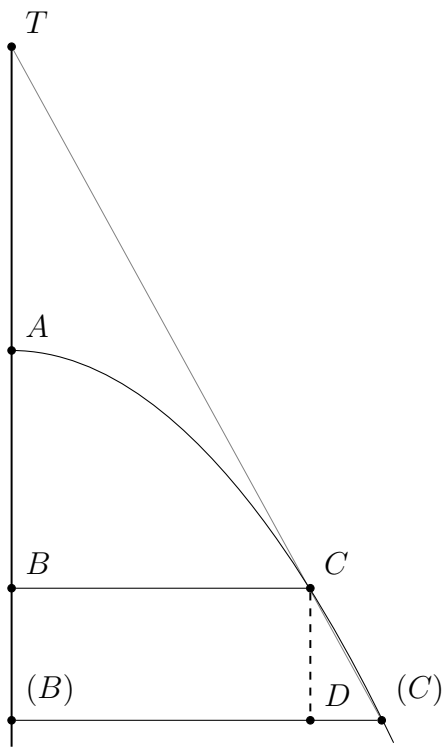

For example, consider Leibniz's analysis, recently translated in
\cite{Le18}, of the inassignable characteristic triangle~$CD(C)$,
where~$D$ is vertex of the right angle whereas~$C$ and~$(C)$ are the
other two vertices; see Figure~\ref{f321}.  This characteristic
triangle is taken to be similar to the assignable triangle~$TBC$.
Leibniz writes:
\begin{quote}
[G]r\^ace \`a ce triangle \emph{inassignable}, c'est-\`a-dire \`a
l'interven\-tion de la raison entre quantit\'es inassignables~$CD$ et
$(C)D$ (que notre calcul diff\'erentiel donne au moyen des quantit\'es
ordinaires ou assignables), on peut trouver la raison entre les
quantit\'es \emph{assignables}~$TB$ et~$BC$, et donc tracer la
tangente~$TC$.  (Leibniz \cite{Le18}, 2018, p.\;155; emphasis added)
\end{quote}
Leibniz's tangent line~$TC$ is thought of as passing through both
infinitely close points~$C$ and~$(C)$.%
\footnote{\label{f25b}To provide a modern interpretation, one can
normalize the equation of the line~$TC(C)$ as~$ax+by=c$
where~$a^2+b^2=1$.  Then one can apply the standard part function (see
note~\ref{f12}) to the coefficients of the equation of the
line~$TC(C)$ to obtain the equation~$a_o x+b_o y=c_o$ of the true
tangent line, where~$r_o=\st(r)$ for each~$r=a,b,c$.  In this sense,
the line~$TC(C)$ and the true tangent line coincide \emph{up to
negligible terms}.  Leibniz often points out that he is working with a
generalized notion of equality.  It the case of the characteristic
triangle, such a notion is applied to secant lines and tangent lines.}
In the same text, Leibniz sees incomparables as equivalent to
inassignables:
\begin{quote}
[O]n voit en optique, quand les divers rayons proviennent d'un m\^eme
point, et que ce point est plac\'e \`a l'infini ou de fa\c con
\emph{inassignable} (ou encore, comme j'ai coutume de dire souvent,
``est \'eloign\'e de fa\c con \emph{incomparable''}), que les rayons
sont parall\`eles.  (ibid.; emphasis added)
\end{quote}
Both involve a violation of Euclid's V.4; see Section~\ref{s23b}.

\subsection{Differentials according to Breger and Spalt}
\label{f20}

Breger writes:
\begin{quote}
I would now like to turn briefly to Leibniz's first publication of his
infinitesimal calculus from 1684. It has been said that Leibniz
introduced infinitesimals here as finite magnitudes (Boyer, 1959, 210;
Bos, 1974, 19, 62--64).  \emph{This is not wrong, but it is
misleading}.  Leibniz in fact explains that one can choose any~$dx$
you like, and he then defines~$dy$ as the magnitude that has the same
relation to~$dx$ as the ordinate to the subtangent.  (Breger
\cite{Br08}, 2008, p.\;188; emphasis added)
\end{quote}
Breger's comments are misleading because they misrepresent Bos's
position.  The first of the pages 62--64 in Bos's article mentioned by
Breger is page\;62.  On this page, Bos is analyzing \emph{Cum
Prodiisset} (\cite{Le01c}, 1701), rather than Leibniz's 1684 article
\cite{Le84} mentioned in Breger's passage.  Here Bos insists on the
difference between the Leibnizian differentials~$dx$ and~$(d)x$, where
the former is inassignable whereas the latter is an ordinary
assignable quantity.  In the passage quoted above, Breger clearly
has~$(d)x$ in mind, but uses the notation~$dx$.  Significantly, Breger
fails to mention anything here about this crucial distinction.
Elsewhere, Breger acknowledges a violation of Euclid V.4 in Leibniz;
see Section~\ref{s34}.

Similarly to Breger, Spalt overlooks the implications of the crucial
Leibnizian distinction between~$dx$ (inassignable) and~$(d)x$
(assignable) in his discussion of Leibnizian differentials when he
writes:
\begin{quote}
Leibniz gives a rigorous \emph{geometric} justification of his rules
for differentials, and, in so doing, to a significant extent he builds
on his law of continuity.  Obviously, Leibniz's differential calculus
has nothing whatsoever to do with a use of infinitely-small
non-variable `numbers', as are known in the modern theory of
nonstandard analysis.  Leibniz' differentials aren't `numbers', but
(variable) \emph{geometrical} `continuous' `magnitudes'.  
(Spalt \cite{Sp15}, 2015, p.\;121; translation ours)
\end{quote}
Spalt raises the issue of the differentials being geometric magnitudes
rather than numbers.  This may be an interesting issue to explore.
However, this issue is transverse to the question of the
non-Archimedean nature of incomparable magnitudes.  Leibniz made it
clear in a letter to l'Hospital that such magnitudes are
non-Archimedean; see Sections~\ref{s22b} and \ref{s34}.

Spalt does mention the Leibnizian differentials~$(d)x$ in a separate
discussion, where he claims that ``The variable length~$dx$ with
limit~$0$ is called `infinitely small'{}'' (Spalt \cite{Sp15}, 2015,
p.\;118).  However, thinking of the Leibnizian~$dx$ as a variable with
limit~$0$ is merely another version of the logical fiction hypothesis.

Spalt goes on to castigate Bos for concluding that 
\begin{quote}
\textsc{Leibniz} had proposed `two distinct' (Bos 1974, p.\;55)
respectively, two `very distinct' (Bos 1980, p.\;70) concepts of
differentials.'' (ibid.)
\end{quote}
Bos' position was outlined in Section~\ref{s11}.  Spalt continues:
\begin{quote}
Advocates of nonstandard analysis routinely refuse to acknowledge
this; the allures of Leibniz' reputation, and of the beautiful field
of activity `historiography of analysis', are too irresistible.%
\footnote{In the original: ``F{\"u}rsprecher der
Nichtstandard-Analysis pflegen sich dieser Einsicht zu
verschlie{\ss}en: Die Verlockungen des Leibniz'schen Renommees sowie
des sch{\"o}nen Bet{\"a}tigungsfeldes ``Geschichtsschreibung der
Analysis'' sind zu unwiderstehlich'' (Spalt\;\cite{Sp15}, 2015,
p.\;121).  }
(ibid.)
\end{quote}
Spalt's move of attributing questionable ideological motives to the
pro-infinitesimal opposition is not without historical precedent.%
\footnote{\label{f21}See e.g., Section~\ref{f13} on Rolle and \emph{la
r\'eforme}.}

\subsection{Infinitesimals and contradictions according to Rabouin}

Rabouin's argument for a version of the logical fiction hypothesis
runs as follows:
\begin{quote}
It should then be clear why infinitesimals were called by Leibniz
``fictions''.  [1]\;In and of itself, there is no such thing as a
``quantity smaller than any other quantity''.  [2]\;This would amount
to the existence of a minimal quantity\; [3]\;and one can show that a
minimal quantity implies \emph{contradiction}\; [4]\;(as a simple
consequence of Archimedes axiom).  [5]\;So ``infinitesimals'' as
``infinitely small quantities'' are terms without reference.  They
only have a contextual meaning and should be paraphrased not by terms,
but by sentences in which only finite quantities occur.  (Rabouin
\cite{Ra15}, 2015, pp.\;362--363; numerals [1] through [5] added;
emphasis added)
\end{quote}
Now it is correct to assert that there is no quantity smaller than any
other quantity, as per Rabouin's sentence~[1].  However, this was not
Leibniz's characterisation of infinitesimals.  Rather, an
infinitesimal is inassignable and is smaller than every
\emph{assignable} quantity (see e.g., Sections~\ref{s1} and
\ref{s24}).  This does not imply the existence of a minimal quantity,
contrary to Rabouin's~[2].  Thus the ``contradiction'' posited by
Rabouin in~[3] is not there.  The reference to the ``Archimedes
axiom'' in [4] is gratuitous; this axiom is not required to obtain a
contradiction starting from Rabouin's incorrect hypotheses.  Thus, his
conclusion~[5] rests on shaky premises.

\subsection{Irrationals, imaginaries, infinitesimals}
\label{s313}

Leibniz employed the qualifier ``impossible'' in reference to
irrational numbers:
\begin{quote}
Des irrationnels naissent les quantit\'es impossibles ou
\emph{imaginaires}, dont la nature est \'etrange, mais dont
l'utilit\'e ne doit pourtant pas \^etre m\'epris\'ee.  En effet,
m\^eme si celles-ci signifient en soi quelque chose d'impossible,
cependant, non seulement elles montrent la source de
l'impossibilit\'e, ainsi que la fa\c con dont la question pourrait
\^etre corrig\'ee afin de ne pas \^etre impossible, mais aussi on
peut, par leur intervention, exprimer des quantit\'es r\'eelles.
(Leibniz \cite{Le18}, 2018, p.\;152; emphasis in the original)
\end{quote}
Leibniz goes on to discuss a few examples, and concludes:
\begin{quote}
Ces expressions ont ceci de merveilleux que dans le calcul elles ne
recouvrent \emph{rien d'absurde ou de contradictoire} et ne peuvent
cependant \^etre montr\'ees dans la \mbox{\emph{nature}},
c'est-\`a-dire dans les choses concr\`etes.  (Leibniz \cite{Le18},
2018, p.\;153; emphasis added)
\end{quote}
Leibniz states clearly that such expressions entail no contradiction.
On the other hand, in ``nature'' there is no referent for such
expressions; cf.\;Bos on infinitesimals as summarized in
Section~\ref{s13}.  Similarly with regard to imaginaries,
\begin{quote}
les valeurs des quantit\'es r\'eelles doivent parfois
n\'ecessaire\-ment \^etre exprim\'ees par l'intervention des
quantit\'es imaginaires et que de l\`a naissent des formules non moins
utiles \`a toute l'\'etendue de l'analyse que ne le sont les formules
communes.  Et ces quantit\'es, je les appelle \emph{impossibles en
apparence}, car \`a la v\'erit\'e elles sont r\'eelles, et je rapporte
les pr\'eceptes par lesquels ceci peut \^etre reconnu.
\cite[pp.\;107--108]{Le18}
\end{quote}
Thus both irrationals and imaginaries are only \emph{apparently}
impossible, according to Leibniz.  Leibniz takes the argument a step
further and makes it clear that such notions imply no contradiction:
\begin{quote}
Il y a une grande diff\'erence entre les quantit\'es imaginaires, ou
impossibles par accident, et celles qui sont absolument impossibles
[parce qu'elles] impliquent contradiction.  \cite[p.\;108]{Le18}
\end{quote}
Leibniz has made it clear that imaginaries (as opposed to
\emph{absolutely impossible} quantities) in fact \emph{do not} imply a
contradiction.  Finally, Leibniz extends the argument to the
infinitely large and small:
\begin{quote}
De fait, les imaginaires ou impossibles par accident, qui ne peuvent
\^etre exhib\'ees parce que fait d\'efaut ce qui est n\'ecessaire et
suffisant pour que se produise une intersection, peuvent \^etre
compar\'ees avec les Quantit\'es infinies et infiniment petites, qui
naissent de la m\^eme fa\c con.  (ibid.)
\end{quote}
Thus according to Leibniz irrationals, imaginaries, and infinitesimals
imply no contradiction; see also Sections~\ref{s25d} and \ref{s45}.
In modern terminology, they are only impossible in the sense of
representing new types of mathematical entities.

\subsection{Apparent impossibility as possibility}
\label{s45}

Leibniz illustrates the \emph{apparent impossibility} of imaginaries
discussed in Section~\ref{s313} via an analysis of a geometric
configuration involving a circle, say~$C$ and a line, say~$L$ (see
\cite{Le18}, p.\;153).  If the nearest distance from~$L$ to the center
of~$C$ is greater than than the radius of~$C$, then~$C$ and~$L$ are
disjoint.  The usual formulas for points in the intersection~$C\cap L$
then contain imaginary terms (in modern language, the intersection
points will appear once the curves~$C$ and~$L$ are complexified).
Leibniz closes the discussion somewhat inconclusively, by commenting
that to create \emph{real} intersection points, one needs to
\emph{change the data of the problem} by either increasing the radius
of~$C$ or moving~$L$ closer to the center of~$C$ (\cite{Le18},
p.\;154).

Leibniz's comment does not solve the problem of accounting for the use
of imaginaries in a situation where one \emph{can't change the data of
the problem}.  Leibniz's comment does not help in situations where one
is forced to make sense of imaginaries and cannot avoid them.  Such a
situation arises e.g., in the solution of the cubic when imaginaries
necessarily arise in an intermediate stage of the calculation, a
technique Leibniz was proficient at; see e.g., Sherry--Katz
(\cite{14c}, 2014, p.\;169).%
\footnote{Here Leibniz is quoted to the effect that ``For this is the
remarkable thing, that, as calculation shows, such an imaginary
quantity is only observed to enter those cubic equations that have no
imaginary root, all their roots being real or possible.''}

Thus Leibniz's discussion of imaginaries in the text translated in
\cite{Le18} (unlike other texts) is incomplete (perhaps it was
completed in subsequent manuscripts), and amounts to walking away from
the problem of the status of imaginaries rather than resolving it.

Leibniz's~$C\cap L$ example is therefore not comparable to his example
of the characteristic triangle (\cite{Le18}, p.\;154--155) discussed
above, where Leibniz does not walk away from the problem but rather
presents a successful solution for finding the tangent line,
\emph{without changing the data of the problem}.  The solution is in
terms of infinitesimals.  Neither example is meant to imply that
imaginaries and/or infinitesimals are either absolutely impossible or
contradictory; on the contrary.  There are many problems treated by
Leibniz where both imaginaries and infinitesimals appear in solutions.

Our analysis undermines Rabouin's claim to the effect that 
\begin{quote}
Le parall\`ele avec les imaginaires est tr\`es souvent mentionn\'e par
ceux qui d\'efendent une vue des infinit\'esimaux comme entit\'es
id\'eales qu'on adjoindrait au domaine des objets r\'eels pour la
r\'esolution de probl\`eme.  Or il est frappant que dans notre texte
comme dans [3b], les imaginaires soient en fait pr\'esent\'es comme
indiquant que le probl\`eme n'a pas de solution.  Si Leibniz pr\'ecise
que les infiniment petits et les points \`a l'infini entrent dans des
probl\`emes qui, eux, ont des solutions, on se gardera n\'eanmoins de
forcer un parall\`ele que Leibniz ne place pas l\`a o\`u
l'interpr\'etation formaliste le place.  (Rabouin in \cite{Le18},
2018, p.\;95, note 2 carrying over to bottom of page 96)
\end{quote}
Contrary to Rabouin's claim, the parallel between imaginaries and
infinitesimals as noncontradictory new types of entities is valid and
is \emph{not} undermined by Leibniz's~$C\cap L$ example, as we
discussed in the current section.

Leibniz repeatedly likens infinitesimals to imaginaries (see also
Sections~\ref{s23} and\;\ref{s22}), and at least once described the
latter as having their \emph{fundamentum in\;re} (basis in fact; see
Leibniz \cite{Le95a}, 1695, p.\;93), providing evidence against the IA
reading that would surely deny them such a basis.

\subsection{Hierarchical structure}

Further support for the reading by Bos is provided by the
\emph{hierarchical} structure on the Leibnizian differentials~$dx$'s,
$dx^2$'s,~$ddx$'s, etc., ubiquitous in Leibniz's texts.  One can
replace~$dx$ by a sequence of finite values~$\epsilon_n$ and furnish a
concealed quantifier incorporated into a modern limit notion so as to
interpret~$dx$ as shorthand for the sequence~$( \epsilon_n)$.
However, one notices that~$\lim_{n\to\infty}\epsilon_n=0$, as well
as~$\lim_{n\to\infty}\epsilon_n^2=0$, and also
unsurprisingly~$\lim_{n\to\infty}(\epsilon_n+\epsilon_n^2)=0$.  Thus,
the Leibnizian substitution
\begin{equation}
\label{e21}
dx+dx^2=dx
\end{equation}
in accordance with his transcendental law of homogeneity (see
Leibniz~\cite{Le10b}, 1710 and also Katz--Sherry \cite{13f}, 2013)
becomes a meaningless tautology~$0+0=0$.  Furthermore, if such
identities are to be interpreted in terms of taking \emph{limits},
then an absurd equality~$dx=dx^2$ would also be true.  To interpret
Leibniz's substitution \eqref{e21} in both a syncategorematic and a
meaningful manner, IA would have to introduce additional ad hoc
proto-Weierstrassian devices with no shadow of a hint in the original
Leibniz.

\subsection{``Historically unforgivable sin''}
\label{s311}

Ishiguro mentions ``Leibniz's followers like Johann Bernoulli, de
l'Hospital, or Euler, who were all brilliant mathematicians rather
than philosophers,'' (Ishiguro \cite{Is90}, 1990, pp.\;79--80) but
then goes on to yank Leibniz right out of his historical context by
claiming that their \emph{modus operandi}
\begin{quote}
is prima facie a strange thing to ascribe to someone who, like
Leibniz, was obsessed with general methodological issues, and with the
logical analysis of all statements and the well-foundedness of all
explanations.  (ibid., p.\;80)
\end{quote}
Having thus abstracted Leibniz from his late 17th--early 18th century
context, Ishiguro proceeds to insert him in a late 19th century
Weierstrassian one.%
\footnote{One of the first occurrences of a modern definition of
continuity in the style of the \emph{Epsilontik} can be found in
Schwarz's summaries of 1861 lectures by Weierstrass; see Dugac
(\cite{Du73}, 1973, p.\;64), Yushkevich (\cite{Yu86}, 1986,
pp.\;74--75).  This definition is a verbal form of a definition
featuring a correct quantifier order (involving alternations of
quantifiers).}
On purely mathematical grounds, such a paraphrase is certainly
possible.%
\footnote{Thus, Robinson noted: ``the method of ultrapowers provides a
ready means for translating a non-standard proof into a standard
mathematical proof in each particular case. However, in the course of
doing so we may complicate the proof considerably, so that frequently
the resulting procedure will be less desirable from a heuristic point
of view. At the same time there may well exist a shorter mathematical
proof which can be obtained independently'' (Robinson \cite{Ro66},
1966, p.\;19).  For a more detailed analysis see Henson--Keisler
(\cite{HK}, 1986).}
However, such an approach to a historical figure would apparently not
escape Unguru's censure:
\begin{quote}
It is \ldots{} a historically unforgiveable sin \ldots{} to assume
wrongly that mathematical equivalence is tantamount to historical
equivalence.  (Unguru \cite{Un76}, 1976, p.\;783)
\end{quote}
Ishiguro seems to have been aware of the problem and at the end of her
Chapter~5 she tries again to explain ``why I believe that Leibniz's
views on the contextual definition of infinitesimals is [sic]
different from those of other mathematicians of his own time who
sought for operationist definitions for certain mathematical notions''
(Ishiguro \cite{Is90}, 1990, p.\;99), but with limited success.
Ishiguro's dubious command of the positions of Leibniz's
contemporaries was discussed in Section~\ref{s23b}.

\subsection{Leibniz against exhaustion}
\label{s312}

Parmentier quotes Leibniz's \emph{De Quadratura Arithmetica} as
follows:
\begin{quote}
J'ai dit jusqu'ici des infinis et des infiniment petits des choses qui
para\^\i tront obscures \`a certains, comme para\^\i t obscure toute
chose nouvelle; rien cependant que chacun ne puisse ais\'ement
comprendre en y consacrant un peu de r\'eflexion pour, l'ayant
compris, en avouer la f\'econdit\'e.  Peu importe que de telles
quantit\'es soient ou non naturelles, on peut se contenter de les
introduire par le biais d'une fiction dans la mesure o\`u elles
offrent bien des commodit\'es
[in the Latin original: \emph{compendia}, meaning ``abbreviations'' or
``shorthand'']
dans les formulations, dans la pens\'ee, et finalement dans
l'invention aussi bien que dans la d\'emonstration, en rendant
inutiles l'usage des figures inscrites et circonscrites, les
raisonnements par l'absurde et la d\'emonstration qu'une erreur est
plus petite que toute erreur assignable.  (Leibniz as translated in
\cite[p.\;284]{Pa01})
\end{quote}
Leibniz's last sentence asserts that infinitesimals make it
unnecessary to get involved in \emph{exhaustion}-type proofs involving
inscribed and circumscribed polygons, arguments with error smaller
than any assignable error, etc.  This would be Leibniz' own refutation
of Ishiguro's thesis.  Namely, Ishiguro claimed that infinitesimals
stand for exhaustion proofs, based on the following passage from
Leibniz:
\begin{quote}
these incomparable magnitudes - are not at all fixed or determined but
can he taken to be as small as we wish in our geometrical reasoning
and so have the effect of the infinitely small in the rigorous sense.
If any opponent tries to contradict this proposition, it follows from
our calculus that the error will be less than any possible assignable
error since it is in our power to make this incomparably small
magnitude small enough for this purpose inasmuch as we can always take
a magnitude as small as we wish.  (Leibniz as translated in
\cite[p.\;87]{Is90})
\end{quote}
However, Leibniz wrote that, on the contrary, infinitesimals make it
unnecessary to get involved in exhaustion proofs.%
\footnote{Note that Ishiguro's sentence ``If magnitudes are
incomparable, they can be neither bigger nor smaller''
\cite[p.\;88]{Is90} involves an equivocation on the term
\emph{incomparable}: if \emph{incomparable} is taken to mean the
definition from the theory of partially ordered sets, then this is a
tautology (roughly ``if magnitudes cannot be compared, then they
cannot be compared''); if \emph{incomparable} is taken to mean ``any
positive-integer multiple is still less than any positive real'' then
Ishiguro's statement is mathematically incorrect, for the hyperreals
are a totally ordered field, hence every element can be compared with
any other.}

\subsection{Leibniz, Des Bosses, and Infinity}
\label{s4b}

Ishiguro quotes a passage from Leibniz's letter to Des Bosses dated 3
march 1706, where Leibniz writes:
\begin{quote}
Meanwhile I have shown that these expressions are of great use for the
abbreviation of thought and thus for discovery as they cannot lead to
error, since it is sufficient to substitute for the infinitely small,
as small a thing as one may wish, so that the error may be less than
any given amount, hence it follows that there can be no error.
(Leibniz as translated in \cite[p.\;85]{Is90})
\end{quote}
Ishiguro concludes: ``It seems then that throughout his working life
as a mathematician Leibniz did not think of founding the calculus in
terms of a \emph{special kind of small magnitude}''
\cite[p.\;86]{Is90} (emphasis added).  But the plain meaning of the
Leibnizian passage is that there are two \emph{distinct} methods, one
involving infinitesimals and one involving errors ``less than any
given amount,'' the former being advantageous over the latter.

Even more significantly, Ishiguro fails to mention that a few months
later, on 1 september 1706, Leibniz wrote another letter to Des Bosses
that sheds light on the question of mathematical infinity.  In this
letter, Leibniz responds to a list of propositions banned by
soon-to-become General Michelangelo Tamburini in 1705.  The list was
sent to Leibniz confidentially by jesuit Des Bosses.  The fourth of
these banned propositions is the following:
\begin{quote}
4. Our minds, to the extent that they are finite, cannot know anything
certain about the infinite; consequently, we should never make it the
object of our discussions.  (translation from Ariew et
al.\;\cite{ACS}, 1998, p.\;258)
\end{quote}
Leibniz comments as follows:
\begin{quote}
Unless I am mistaken, mathematicians have already refuted the fourth
proposition, and I myself have published some samples of the science
of the infinite.  However, I maintain, strictly speaking, that an
infinite composed from parts is neither one nor a whole, and it is not
conceived as a quantity except through a fiction of the mind.
(Leibniz as translated in \cite{Le07}, 2007, p.\;53)
\end{quote}
Here Leibniz affirms that human mind can indeed conceive of infinity
(contrary to proposition 4 rejected by both the jesuits and himself),
and moreover that he published ``samples of the science of the
infinite" to prove this.  Here clearly mathematical infinity is not a
mere sign for hidden quantifiers, contrary to Ishiguro's position.

The preponderance of the evidence in the primary sources indicates
that Leibniz did indeed found his calculus on a special kind of
fictional small magnitude.

\section{\emph{De Quadratura Arithmetica}}
\label{s4}

Leibniz's unpublished text \emph{De Quadratura Arithmetica}\ldots (DQA)
was written shortly after he developed the infinitesimal calculus in
1675.  Thus the work dates from an early period of his mathematical
career.  Here Leibniz wrote:
\begin{quote}
Nec refert an tales quantitates sint in rerum natura, sufficit enim
fictione introduci, cum loquendi cogitandique, ac proinde inveniendi
pariter ac demonstrandi compendia praebeant, ne semper inscriptis vel
circumscriptis uti\ldots necesse sit.%
\footnote{Jesseph's translation: `` `Nor does it matter whether there
are such quantities in nature, for it suffices that they are
introduced as fictions, since they allow the abbreviations of speech
and thought in the discovery as well as demonstration' (Leibniz 1993,
p.\;69)'' (Jesseph \cite{Je15}, 2015, p.\;198).  A longer passage
including this one was quoted in Parmentier's French translation in
Section~\ref{s312}.}
(Leibniz \cite{Le93}, p.\;69)
\end{quote}

\subsection{B-track reading}
\label{s41}

A straightforward interpretation of this passage from DQA is that
there exist two approaches to the calculus:
\begin{enumerate}
\item[(A)] one involving inscribed and circumscribed figures, called
the method of exhaustion; and
\item[(B)] a method involving what he referred to elsewhere as
\emph{useful fictions}, and enabling abbreviations of speech and
thought when compared to method\;A.
\end{enumerate}
The theme of a pair of distinct approaches occurs often in Leibniz's
writing.  Thus, in his 2\;february\;1702 letter to Varignon, Leibniz
writes:
\begin{quote}
Et c'est pour cet effect que j'ay donn\'e un jour des lemmes des
incomparables dans les Actes de Leipzic, qu'on peut entendre comme on
vent%
\footnote{I.e., ``veux''.}
[sic], soit des infinis \emph{\`a la rigueur}, soit des grandeurs
seulement, qui n'entrent point en ligne de compte les unes au prix des
autres.  Mais il faut considerer en m\^eme temps, que ces
incomparables communs m\^emes n'estant nullement fixes ou
determin\'es, et pouvant estre pris aussi petits qu'on veut dans nos
raisonnemens Geometriques, font l'effect des infiniment petits
rigoureux\ldots{} (Leibniz \cite{Le02}, 1702, p.\;92; emphasis added)
\end{quote}
The passage is analyzed in its context in Section~\ref{s25c}.

\subsection{A-track reading}

An alternative interpretation following Ishiguro of this type of
passage in Leibniz is that the B-method is merely shorthand for the
A-method involving hidden quantifiers, in the spirit of Russell's
\emph{logical fictions}.  Such an interpretation of the passage seems
more forced than the one we gave above.%
\footnote{\label{f24}Based on the \emph{second} sentence of the
Leibnizian passage quoted in Section~\ref{s41}, Breger claims that
``[Leibniz] stresses that the incomparably small magnitudes are
certainly not `fixes ou determin\'es'; they can be chosen as small as
one wants'' (Breger \cite{Br17}, 2017, p.\;77).  However, Breger fails
to mention the fact that the \emph{first} setence of the passage we
quoted makes it clear that Leibniz's comments on incomparables not
being ``fixes ou determin\'es'' apply to \emph{common} incomparables
of track A, rather than the ``infinies \`a la rigueur'' of track B.}

In DQA Leibniz wrote that infinitely small (say,~$\epsilon$) means
smaller than any given quantity, and that infinitely large means
larger than any given quantity.

Such descriptions similarly admit two readings: (A) for every given
positive quantity there is an~$\epsilon>0$ smaller than it; or (B)
this specific~$\epsilon$ is smaller than every given (i.e.,
assignable) positive quantity.  In reading~(B) the~$\epsilon$ is
\emph{inassignable}.  The dichotomy of assignable \emph{vs}
inassignable quantity (or magnitude) was used frequently by Leibniz,
e.g., in \emph{Cum Prodiisset} \cite{Le01c} and \emph{Puisque des
personnes\ldots} \cite{Le05b}%
\footnote{See also note~\ref{f4} and Section~\ref{s24} on the
dichotomy assignable \emph{vs} inassignable.}

Leibniz's reference to Archimedes in various texts typically is a
reference to the method of exhaustion, and is sometimes accompanied by
a claim that infinitesimals violate Euclid Definition~V.4 (see e.g.,
Section~\ref{s22b}).  The latter is called today the Archimedean
property (but was not in Leibniz's time).

The term \emph{syncategorematic} itself is incidental to the true
issues involved.  The real issue is whether Leibniz was relying on
hidden quantifiers (as per IA) or not.

\subsection
{Bl\aa sj\"o--Knobloch differences over DQA}

Knobloch wrote in reference to DQA:
\begin{quote}
In this treatise, Leibniz laid the rigorous foundation of the theory
of infinitely small and infinite quantities or, in other words, of the
theory of quantified indivisibles. In modern terms Leibniz introduced
`Riemannian sums'%
\footnote{The correct technical term for this concept is \emph{Riemann
sum}.  The adjective \emph{Riemannian} is also used in a technical
sense, but in other contexts; e.g., \emph{Riemannian geometry}.}
in order to demonstrate the integrability of continuous functions.
(Knobloch \cite{Kn02}, 2002, p.\;59)
\end{quote}
Knobloch then proceeds to describe Leibniz's method, and notes:
\begin{quote}
While the ``common method of indivisibles''%
\footnote{\label{f25}See note~\ref{f2}.}
considered inscriptions and circumscriptions of mixtilinear figures,
the step figure is neither an inscription nor a circumscription,
rather something in between. In modern terms: Leibniz demonstrated the
integrability of a huge class of functions by means of Riemannian sums
which depend on intermediate values of the partial integration
intervals.  (ibid., p.\;63)
\end{quote}
Thus, Knobloch argues that Leibniz's technique in DQA represented an
advance over earlier inscription and circumscription techniques.

Jesseph argues in \cite{Je15} that the techniques in DQA were limited
by their reliance on the knowledge of the tangent lines to the curve
(and/or the corresponding differential).  Therefore the applicability
of the techniques depended on the availability of such data.%
\footnote{Thus, Leibniz applies his method in DQA to find the
quadrature of general cycloidal segments (Edwards \cite{Ed79}, 1979,
p.\;251).  Here the calculation exploits the
%
%
geometric knowledge related to the construction of the curve.}
Accordingly, the solution of the quadrature problem in DQA depends on
a differentiation problem.

Bl\aa sj\"o analyzes Leibniz's technique in DQA as follows.  The
so-called general integration theorem in DQA assumes the existence of
tangent lines, not only of the ``general function'' one starts with,
say~$f(x)$, but also for a secondary function, say~$d(x)$ (following
the notation in Bl\aa sj\"o \cite{Bl17a}, 2017), is constructed
from~$f(x)$ by means of the tangents to~$f(x)$.  The theorem relates
the areas under~$f$ and~$d$ to each other in a manner closely related
to modern integration by parts.  The assumption that~$f(x)$ has
tangents everywhere (or possibly almost everywhere) is essential since
otherwise there wouldn't even be any function~$d(x)$ to investigate.
The assumption that~$d(x)$ has a tangent everywhere is less essential.
The \emph{knowledge} of the tangents themselves does not play an
essential role in the proof, but it is essential that the curve has no
``reversion points", which is a notion that Leibniz has not defined
otherwise than in terms of tangents (and hence assumed that tangents
must exist, for a non-differentiable curve could reverse directions
without having a tangent corresponding to the turning point).  This
concerns the general condition under which the theorem is valid.  Was
Leibniz explicating precise and rigorous conditions of validity?
Clearly he was not.  The conditions he does state are of an intuitive
nature and are not intended as rigorous conditions of validity.  For
the latter purpose, they are clearly insufficient, for example because
of the issue of the tangents.  Leibniz's result does not provide any
general integration theory, and is actually a single specific
integration technique, closely related to integration by parts.  Just
as the modern integration by parts formula, it can be considered to
apply in great generality, but it is only useful in cases where we
have some
%
%
geometric information equivalent to the derivatives or antiderivatives
involved (otherwise it is just an exercise in expressing some unknown
integral in terms of another unknown integral).  For more details see
Bl\aa sj\"o (\cite{Bl17a}, 2017).

Our thesis in the present text is independent of resolving these
differences among scholars concerning DQA.  Namely, the contention
that Leibniz's infinitesimal was inassignable in the sense of
violating the Archimedean property when compared to ordinary
(assignable) quantities is independent of the Bl\aa sj\"o--Knobloch
differences over DQA.  Whatever the foundational significance of DQA
may have been (and this is the subject of their disagreement), the
fact remains that here Leibniz talks about two separate methods: track
A and track B.  Leibniz seeks to justify the direct track-B method
(exploiting inassignable infinitesimals) in terms of an
exhaustion-type track-A method.

\subsection{Logarithmic curve}

Leibniz also investigated the special case of the logarithmic curve in
Proposition 46 in DQA.  In the statement of Proposition 46, Leibniz
speaks of information that can be obtained from the hyperbola in terms
of which the logarithmic curve was defined.  Thus the investigation
depends on the knowledge that~$\log(x)$ is the area
under~$\frac{1}{x}$ (in other words, the knowledge of the derivative
of the logarithm).  This can be easily understood in the context of
integration by parts.  Therefore Leibniz's propostion on the
logarithmic curve is not an example of performing quadratures without
knowledge of derivatives or tangents.

\subsection{Manuscript remained unpublished}
Leibniz never published the DQA.%
\footnote{The loss of a manuscript version in transit from Paris to
Hannover in 1679, signaled by Knobloch \cite[p.\;282]{Kn17}, could
have been overcome by writing a new version, Leibniz having been only
33 at the time.}
Was that because he realized that the A-method, while not
contradicting the B-method, was an impediment to the \emph{Ars
Inveniendi}?  Jesseph notes in \cite[p.\;200]{Je15} that Leibniz may
have set aside the DQA without publishing it because he had turned his
attention to more powerful methods that he would introduce in the
1680s in what he called ``our new calculus of differences and sums,
which involves the consideration of the infinite", and ``extends
beyond what the imagination can attain"%
\footnote{The full sentence in the original reads: ``Ainsi il ne faut
point s'\'etonner, si notre nouveau calcul des diff\'erences et des
sommes, qui enveloppe la consid\'eration de l'infini et s'\'eloigne
par cons\'equent de ce que l'imagination peut atteindre, n'est pas
venu d'abord \`a sa perfection.''  See Gerhardt's edition
(\cite{Ge50}, vol.\;V, p.\;307).  The passage appears in Leibniz's
article ``Considerations sur la diff\'erence qu'il y a entre l'analyse
ordinaire et le nouveau calcul des transcendantes'' in \emph{Journal
des S\c cavans} in 1694.}
(GM V, 307).

\section{Conclusion}
\label{s5}

Leibniz did not merely use an infinitesimal approach as an unrigorous
way of doing calculus that makes the work easier.  Rather, what we are
arguing is that historically there have been two separate approaches
to the calculus: track (A) and track (B).  The historical calculus has
often been criticized from an anachronistic modern set-theoretic
viewpoint that would make both approaches appear unrigorous to the
extent that they did not possess a set-theoretic justification which
is considered a \emph{sine-qua-non} of rigor in today's mathematics.
From such an anachronistic standpoint, \emph{both} historical
approaches were unrigorous by today's standards.

The work of Fermat, Gregory, Leibniz, Euler, Cauchy and others created
a body of \emph{procedures} called infinitesimal calculus and/or
analysis, in what could be referred to as the pioneering phase of the
discipline.  Following the pioneering phase, efforts were made to
develop set-theoretic justification for this body of procedures.
Eventually this effort succeeded both for track (A) and track (B).
Assigning names is a matter of debate but in a sense Edward Nelson's
syntactic approach (\cite{Ne77}, 1977) is particularly fundamental
because it shows that one can take an infinitesimal to be a primitive
notion within the context of the ordinary real line, in the spirit of
what natural philosophers since at least Pascal have envisioned.
Nelson's approach exploits a unary (i.e., one-place) predicate
\textbf{standard}; the formula~$\mathbf{standard}(x)$ reads ``$x$ is
standard''.  Thus, mathematical entities can be either \emph{standard}
or \emph{nonstandard}.  This applies in particular to real numbers.
The standard/nonstandard distinction can be seen as a formalisation of
Leibniz's assignable/inassignable distinction.%
\footnote{\label{s32}In Nelson's approach, the violation of the
Archimedean property takes the form~$(\exists \varepsilon>0)
(\forall^{\st}n\in\N)\, \left[\varepsilon <\frac{1}{n} \right]$, where
$\forall^{\st}$ is universal quantification over \emph{standard}
elements only.}
For details on the systems developed by Nelson, Hrbacek, Kanovei, and
others see Fletcher et al.\;(\cite{17f}, 2017).

Ishiguro and Arthur have argued that what appears to be a B-track
method is in reality a Russellian illusion that is eliminable by
careful analysis of the pioneering texts of Leibniz.  This is the
meaning of their \emph{logical fiction} hypothesis; Arthur goes so far
as to speak, in the title of his forthcoming chapter (\cite{Ar19},
2019), of ``Archimedean infinitesimals,'' a word string that would
appear as a freshman quantifier-order error to many a mathematically
educated scholar.  We hold the IA hypothesis to be an error of
interpretation and have argued that it is not backed by solid textual
evidence in Leibniz.

Leibniz, by arguing in favor of exploiting \emph{inassignable}
infinitesimals even though they are fictions, differed from l'Hopital
and Bernoulli who were prepared to assign a loftier ontological status
to infinitesimals as truly existing entities.  In this Leibniz was
remarkably modern, and anticipated formalist strategies that fully
emerged in the 20th century.  As C. H. Edwards, Jr. points out,
\begin{quote}
It is important to note that the differentials~$dx$ are fixed
non-zero quantities; they are neither variables approaching zero nor
ones that are intended to eventually approach zero.
\cite[p.\;261]{Ed79}
\end{quote}
Many mathematicians, historians, and philosophers nowadays are in
favor of pluralism as far as foundations of mathematics are concerned.
Today we have both A-type and B-type set-theoretic foundations, as
well as category-theoretic foundations both of classical and
intuitionistic types.  The historical studies currently available
suggest that in the case of Leibnizian infinitesimal calculus, A-type
foundations are insufficiently expressive to capture the spirit of
Leibniz's work.  Pluralism is a good thing in principle but the A-type
logical fiction interpretation is not a viable alternative to the
B-type pure fictional one.

The \emph{logical fiction} reading of Leibnizian infinitesimals has
become entrenched to an extent that some Leibniz scholars feel
compelled to endorse it publicly, while in private correspondence
conceding that Leibniz used a dual strategy which we have elaborated
in terms of a distinction between Leibnizian methods (A) and (B) (see
Section~\ref{s41}), the latter involving \emph{pure fictional}
infinitesimals, as opposed to the IA \emph{logical fiction}
hypothesis.  The Ishiguro--Arthur hypothesis must be rejected as
having little basis in Leibniz's writings.

\section*{Acknowledgments} 

We are grateful to John Dawson for providing a copy of Robinson's
23\;august\;1973 letter to G\"odel quoted in Section~\ref{s37}; to
Roger Ariew for bringing to our attention Leibniz's letter to Des
Bosses dated 1\;september 1706, analyzed in Section~\ref{s4b}; to
Charlotte Wahl for bringing to our attention Leibniz's manuscript
\emph{Puisque des personnes\ldots} (\cite{Le05b}, 1705), analyzed in
Section~\ref{s1}; to Vladimir Kanovei, Eberhard Knobloch, Semen
Kutateladze, Marc Parmentier, Siegmund Probst, and David Rabouin for
helpful comments on earlier versions of the manuscript; and to Viktor
Bl\aa sj\"o for helpful comments on Leibniz's DQA.\, M.\;Katz was
partially supported by the Israel Science Foundation grant
no.\;1517/12.

\end{document}